\def\ZZ         {{\bf Z}}
\def\RR         {{\bf R}}
\def\CC         {{\bf C}}
\def\QQ         {{\bf Q}}
\def\mm         {{\bf m}}
\def\nn         {{\bf n}}

\documentstyle[twoside,12pt]{article}

\newtheorem{prop}{\rm PROPOSITION}[section]
\newtheorem{dfn}[prop]{\rm DEFINITION}

\newtheorem{rem}[prop]{\rm REMARK}

\newtheorem{lem}[prop]{\rm LEMMA}

\title{A FINITENESS THEOREM FOR SUBGROUPS OF SP(4,\ZZ)}

\author{
Lev A. Borisov
\\
\small Department of Mathematics,  Columbia University \\
\small 2990 Broadway, Mailcode 4432, New York, NY 10027, USA \\
\small e-mail: lborisov@math.columbia.edu}

\begin{document}

\date{}

\maketitle

\begin{abstract}
This paper proves that there are only finitely many subgroups $H$ of 
finite index in ${\rm Sp(4,\ZZ)}$ such that the corresponding quotient 
${\cal H}/H$ of the Siegel upper half space of rank two is not of general
type.
\end{abstract}

\section{Introduction}

The Siegel upper half space of rank two consists of complex symmetric two by two
matrices whose imaginary part is positive definite. It will be denoted by
${\cal H}$ throughout the paper. It is the moduli space of principally
polarized marked abelian surfaces. 
The group ${\rm Sp(4,\ZZ)}$ acts on ${\cal H}$ by the automorphisms of the 
marking. This group consists of four by four integer matrices of the form
$\pmatrix{A&B \cr C&D\cr}$ where $A,B,C,$ and $D$ are two by two matrices that 
obey $A{}^tB=B{}^tA,C{}^tD=D{}^tC,A{}^tD-B{}^tC={\bf 1}.$
Written in coordinates, this action becomes
$$\pmatrix{A&B\cr C&D\cr}\cdot M=(AM+B)(CM+D)^{-1}.$$
It is a natural generalization of the usual upper half plane with the action
of ${\rm Sl(2,\ZZ)}$. It is related to various moduli spaces of 
abelian surfaces in the same way the usual upper half plane is related 
to moduli spaces of elliptic curves.

We shall be concerned mostly with quotients of ${\cal H}$ by 
the action of subgroups $H$ of finite index in ${\rm Sp(4,\ZZ)}.$
These quotients are known to be algebraic varieties of dimension $3$.
They have been studied extensively since the end of last century.
Some of these varieties have extremely rich
and beautiful geometry, see for instance \cite{Geer},\cite{Lee} and \cite{GeerII}.

The goal of this paper is to prove the following statement, see
proposition \ref{fintheorem}.

{\bf Finiteness Theorem.} There are only finitely many subgroups 
$H\subseteq {\rm Sp(4,\ZZ)}$ of finite index such that ${\cal H}/H$
is not of general type.

The important corollary of this result is that there are only finitely many
subgroups $H$ such that the quotient ${\cal H}/H$ is rational.
Varieties of general type can be viewed as the generalization to higher 
dimension of curves of genus two or more.  
It is reasonable to expect that they do not have any special geometric properties, 
and thus all interesting quotients ${\cal H}/H$ can be in 
principle listed. This theorem is analogous to the result of J.G. Thompson 
(see \cite{Thompson}) for the usual upper half plane.
More accurate estimates for certain classes of subgroups of ${\rm Sp(4,\ZZ)}$
have been proved in \cite{Grady,Gritsenko,Hulek}.

The method of the proof is roughly the following. It is known that $H$
contains a principal congruence subgroup $\Gamma(n)$ of some level $n$.
The quotient ${\cal H}/\Gamma(n)$ admits a well understood smooth 
compactification, constructed in the paper of Igusa \cite{Igusa}.
Our aim is to construct global sections of the multicanonical
line bundle on the desingularization of the compactification of 
$({\cal H}/\Gamma(n))/(H/\Gamma(n))$
from the sections of certain line bundles on the Igusa compactification of
${\cal H}/\Gamma(n)$. 

We will use standard facts about singular algebraic varieties, which are
collected in Section 7. The results of Sections 2 and 4 
are probably known to specialists in the field, although there
are not many convenient references. Section 3 and 5 are the key
sections of the paper. The former is a purely combinatorial calculation, and
the latter is an algebra-geometrical one. In both sections we assume
that $n$ is a power of a prime, and Section 6 allows
us to drop this restriction.

This paper is essentially my University of Michigan thesis. 
Major part of it was done when I was still in Moscow. 
It is influenced a lot by my advisor Vasilii Iskovskikh
who taught me the basics of algebraic geometry as well as some 
singularity theory which comes in very handy in the paper. 
I would like to thank 
Osip Shvartsman for many stimulating discussions on the subject of this paper.
My thesis advisor Igor Dolgachev has been a constant source of inspiration for my
studies of algebraic geometry at the University of Michigan. I also wish
to thank
Gopal Prasad for several valuable conversations and Melvin Hochster for providing
a useful reference.

\section{Algebraic cycles on Satake and Igusa compactifications}

The purpose of this section is to recall the basic facts about some 
special algebraic cycles on the Satake and Igusa compactifications of
${\cal H}/\Gamma(n)$ and to find a nice combinatorial description of
their components.

We consider the principal congruence subgroup $\Gamma(n)$ of level
$n$ inside ${\rm Sp(4,\ZZ)}$. For the rest of the section $n$ is fixed
and is greater than two. The group $\Gamma(n)$ acts on the Siegel upper
half space of rank two ${\cal H}$ according to the formula
$$\pmatrix{A&B\cr C&D\cr}\cdot\tau=(A\tau+B)(C\tau+D)^{-1}.$$

The quotient ${\cal H}/\Gamma(n)$  is a nonsingular algebraic variety. 
It is a Zariski open subset of the compact singular algebraic 
variety called the Satake compactification of ${\cal H}/\Gamma(n)$.
The exact references can be found in \cite{Igusa}. The monoidal 
transformation of the Satake compactification along the singular
locus is nonsingular. This
variety was first considered by Igusa in \cite{Igusa}, and is called
the Igusa compactification. We denote it by $X_n$. Points of
${\cal H}/\Gamma(n)$ are referred to as the {\it finite}\/ part of the 
compactification and the rest is the part {\it at infinity.}

The part at infinity of the Satake compactification consists of a finite
number of curves that intersect in a finite number of {\it cusp}\/
points. The part at infinity of the Igusa compactification is a divisor
$D=\sum_iD_i$, which has simple normal crossings. Its components are elliptic
fibrations over the curves at infinity of the Satake compactification.
The group $G=\Gamma(1)/\Gamma(n)$ acts on both compactifications,
and the map between them is equivariant. The group $G$ is isomorphic
to ${\rm Sp(4,\ZZ/n\ZZ)}$, and $\pm{\bf 1}$ act as the identity.

There are two more types of divisors on the Igusa compactification 
that will be important to our discussion. First of all, there are
divisors $E_i$ that are conjugates of the closure of the image
in ${\cal H}/\Gamma(n)$ of the set of diagonal matrices in 
${\cal H}$. They are disjoint and are isomorphic to the product of
two modular curves (see \cite{Yamazaki}). We denote their sum by $E$.
We also consider divisors that are conjugates of
the closure of the image of the set of matrices $\pmatrix{x&y\cr y&x\cr}$
in ${\cal H}$. Geometrically, these matrices correspond to Jacobians
of genus two curves with an extra involution, see \cite{Bolza}.
We denote them  by $F_i$ and their sum by $F$. They do intersect with 
each other and their geometry is somewhat more complicated. We prove
the necessary statements regarding these at the end of this section.

We abuse notation somewhat to denote ${\rm Sp(4,\ZZ)}$-conjugates
of the sets $\pmatrix{x&0\cr 0&z\cr}$ and $\pmatrix{x&y\cr y&x\cr}$ by
$E_i$ and $F_j$ as well.

Let us introduce the abelian group $V$ of column vectors of
height four with coefficients in $\ZZ/n\ZZ$
provided with the skew form $\langle ~,~\rangle $ defined by the formula
$\langle {}^t(x^1,...,x^4),{}^t(y^1,...,y^4)\rangle =x^1y^3+x^2y^4-x^3y^1-x^4y^2.$
The group $G$ acts naturally on $V$ by left multiplication.
Our goal here is to construct $G$-equivariant
correspondences between components of cycles on the Satake
and Igusa compactifications mentioned above and some objects defined 
in terms of the group $V$.

\begin{prop} 
{ The infinity divisors of the Igusa compactification (or equivalently,
the curves at infinity of the Satake compactification) are in one-to-one $G$-
equivariant correspondence with the primitive $\pm$vectors $\pm v$ in $V$.
Here we call a vector $v$ primitive iff its order is exactly $n$. 
The $\pm$ means that we identify opposite vectors.}
\label{indexD}
\end{prop}

{\em Proof.} It is known (see \cite{Igusa}) that all components of $D$ are
$G$-conjugate. It can be shown that the group $G$ also acts transitively on 
the set of primitive $\pm$vectors. It remains to notice that the stabilizer
of the $\pm$vector ${}^t(0,1,0,0)$ coincides with the stabilizer of $D_0$,
where $D_0$ is the {\it standard}\/ divisor that corresponds to the basis
of open subsets $\{\pmatrix{x&y\cr y&z\cr},~{\rm Im}(z)\to+\infty\}$ of ${\cal H}.$
The description of the stabilizer of $D_0$ can be derived from \cite{Igusa}.
It consists of matrices of the form 
$$\pm\pmatrix{a&0&b&m_3\cr m_1&1&m_2&m_4\cr
 c&0&d&m_5\cr 0&0&0&1\cr}({\rm mod}n),$$$$~ad-bc=1({\rm mod}n),
 ~bm_1+m_3=am_2({\rm mod}n),~dm_1+m_5=cm_2({\rm mod}n).$$

This allows us to construct a bijective correspondence between infinity
divisors on the Igusa compactification and $\pm$vectors in $V.$
We will use the notation $\pm v_\alpha$ for the $\pm$vector that corresponds
to the divisor $D_\alpha$ and vice versa.\hfill$\Box$

\begin{prop}
{ Cusp points $Q_i$ of the Satake compactification
are in one-to-one $G$-equivariant correspondence with the following pairs
 $(W,\pm f)$. We consider all possible $W\subset V$ and $f:W\times W\to 
{\ZZ/n\ZZ}$ such that

(1) $W$ is a subgroup of $\,V$ isomorphic to $({\ZZ/n\ZZ)}^2$,

(2) $\langle ,\rangle |_W=0$,

(3) $f$ is a non-degenerate skew form on $W$ with values in 
${\ZZ/n\ZZ}$, where non-degeneracy means $f(W\times W)\ni 1(n)$.} 
\label{indexptsSatake}
\end{prop}

{\em Proof.} All cusp points are conjugates of the one described 
by the basis of open sets  $$\{\pmatrix{x&y\cr y&z\cr},~{\rm Im}(\pmatrix{x&y\cr
y&z\cr})\to+\infty\}$$
(see \cite{Igusa}). We call this point {\it standard}\/. 
The stabilizer of the standard point consists of matrices
of the form
$$\{\pmatrix{A&B\cr {\bf 0}&{}^tA^{-1}\cr}, A{}^tB=B{}^tA, det(A)=\pm1(n)\}$$
However, this is exactly the stabilizer of the {\it standard}\/
pair $$(W,\pm f)=({}^t(*,*,0,0),f({}^t(1,0,0,0),{}^t(0,1,0,0))=1(n)).$$
It can be shown that any pair $(W,f)$ is a $G$-conjugate of the standard pair.
As a result, we can define the required $G$-equivariant correspondence.
We will use the notation $(W_\alpha,\pm f_\alpha)$ for the pair that 
corresponds to the point $Q_\alpha$ and vice versa. \hfill$\Box$

\begin{prop} 
{ The curve at infinity of the Satake compactification that corresponds
to the divisor $D_\alpha$ contains the cusp point $Q_\beta$ iff
$v_\alpha\in W_\beta$.}
\label{indexcurvethroughpointonSatake}
\end{prop}

{\em Proof.} Consider the action of the group that stabilizes the standard
curve. It acts transitively on the set of cusp points of this curve, which
are exactly the $Q_i$'s. Therefore, all inclusion pairs are acted upon
transitively. The standard curve passes through the standard point,
and ${}^t(0,1,0,0)\in {{}^t(}*,*,0,0)$, which proves the only if part of the
statement. On the other hand, the stabilizer of the standard point acts
transitively on the $\pm$vectors in ${}^t(*,*,0,0)$, which proves the if
part. \hfill$\Box$

\begin{prop}
{ Two infinity divisors $D_\alpha$ and $D_\beta$ intersect over
the point $Q_\delta$ iff $v_\alpha,v_\beta \in W_\delta$ and
$f_\delta(v_\alpha,v_\beta)=\pm 1(n)$. In this case the intersection is
isomorphic to ${\bf P}^1$.}
\label{indexDD}
\end{prop}

{\em Proof.} Because of transitivity of the action, the point $Q_\delta$
may be considered standard. We follow the argument of \cite{Igusa} for the 
case where $g_0=0$ and $g_1=2$. Curves of the intersection
of the two infinity divisors are conjugate to one of the curves obtained by
taking the limits of the points $\pmatrix{x&y\cr y&z\cr}$,
with imaginary parts of two out of three normal coordinates
$y,(-x-y),(-z-y)$ going to $-\infty$ and the remaining one being bounded.
They are pairwise intersections of the divisors that correspond to
the limits where exactly one of the imaginary parts goes to $-\infty$
and the other two are bounded. The divisor that corresponds to
$Im(z+y)\to\infty$ is exactly the standard divisor, because ${\rm Im}(y)$
is bounded. The other two divisors are obtained from the standard one
by the action of $\{\pmatrix{A&{\bf 0}\cr {\bf 0}&{}^tA^{-1}\cr}\}$ with
$A=\pmatrix{1&-1\cr 0&1\cr}$ and $A=\pmatrix{0&1\cr 1&0\cr}$
respectively. Therefore, these divisors correspond to the $\pm$vectors
$\pm {}^t(0,1,0,0)$, $\pm {}^t(-1,1,0,0)$, and $\pm {}^t(1,0,0,0).$
This proves the "only if" part. The "if" part follows from the transitivity
of the action of $G$ on the combinatorial data on the right hand side of
the statement. The fact that each irreducible component of the intersection
is isomorphic to ${\bf P}^1$ is proven in \cite{Igusa}, and the uniqueness of
the irreducible component can be derived easily from the description of
the divisors in terms of the above limits.
\hfill$\Box$

\begin{prop} 
{ Three infinity divisors $D_\alpha, D_\beta, D_\gamma$ intersect over
the point $Q_\delta$ iff $v_\alpha,v_\beta,v_\gamma \in W_\delta$,
the set $\{\pm v_\alpha \pm v_\beta \pm v_\gamma\}$ contains $0$, and
$f_\delta(v_\alpha,v_\beta)=\pm 1(n)$. In this case the intersection point
is unique.}
\label{indextDDD}
\end{prop}

{\em Proof.} As in the previous proposition, we prove that all points
of triple intersection are conjugates of the intersection point of the 
divisors that correspond to $\pm {}^t(0,1,0,0)$, $\pm {}^t(-1,1,0,0)$, 
and $\pm {}^t(1,0,0,0).$ Then we notice that any triple of $\pm$vectors
with the above properties can be transformed to the triple
$(\pm {}^t(0,1,0,0)$, $\pm {}^t(-1,1,0,0)$, $\pm {}^t(1,0,0,0)).$
\hfill$\Box$

\begin{prop}
{ Divisors $E_i$ are in one-to-one $G$-equivariant correspondence
with unordered pairs $(W_1,W_2)$ such that 

(1) $W_1$ and $W_2$ are subgroups of $V$ isomorphic to 
$({\ZZ/n\ZZ})^2$ each,

(2) $W_1 \perp W_2=V$.}
\label{indexE}
\end{prop}

{\em Proof.} All divisors $E_i$ are conjugates of the {\it standard}\/
one defined as the closure of the image of the set of diagonal matrices. 
The stabilizer of this standard divisor is described in \cite{Yamazaki}.
It is exactly the stabilizer of the {\it standard}\/ pair
$({}^t(*,0,*,0),{}^t(0,*,0,*)).$ It is easy to show that every pair $(W_1,W_2)$
with above properties is conjugate to this standard one, which
completes the proof. For a given $E_\alpha$ the corresponding pair
will be denoted by $(W_{\alpha1},W_{\alpha2})$ and vice versa.
\hfill$\Box$

\begin{prop}
{ The divisor $E_\alpha$ intersects the divisor $D_\beta$ iff
$v_\beta$ lies in one of the subgroups $W_{\alpha1},W_{\alpha2}$.
In this case the intersection is isomorphic to the modular curve
of principal level $n$.}
\label{indexED}
\end{prop}

{\em Proof.} We assume that the divisor $E_\alpha$ is a standard one.
Then the statement of the proposition follows from the description of
the action of the group $\Gamma(n)$ in a neighborhood of the set of
diagonal matrices (see \cite{Yamazaki}).
\hfill$\Box$

There is an alternative way to describe divisors $E_i$.

\begin{prop} 
{ Divisors $E_i$ are in one-to-one $G$-equivariant correspondence
with conjugates of the involution
$$\pm\pmatrix{1&0&0&0\cr 0&-1&0&0\cr 0&0&1&0\cr 0&0&0&-1\cr}$$
in the group ${\rm Sp(4,\ZZ/n\ZZ)/\{\pm 1\}}.$}
\label{indexinvE}
\end{prop}

{\em Proof.} The action of 
$$\pm\pmatrix{1&0&0&0\cr 0&-1&0&0\cr 0&0&1&0\cr 0&0&0&-1\cr}$$
on ${\cal H}$ is defined by $\pmatrix{x&y\cr y&z\cr}\to \pmatrix{x&-y\cr 
-y&z\cr}$, so it fixes exactly the points of the standard divisor
$E_0.$ This gives a one-to-one correspondence between $\Gamma(1)$
conjugates of this involution and $\Gamma(1)$ conjugates of the 
diagonal in ${\cal H}.$ This correspondence survives when we mod out
by $\Gamma(n)$, and then we use surjectivity of 
$\Gamma(1)/\Gamma(n)\to {\rm Sp(4,\ZZ/n\ZZ)}.$
\hfill$\Box$

This alternative description is related to the original one as follows.

\begin{prop}
{ The involution $\varphi_\alpha$ that fixes all points of the
divisor $E_\alpha$ is defined by

(1) $\varphi_\alpha |_{W_{\alpha1}}={\rm id} |_{W_{\alpha1}}$,

(2) $\varphi_\alpha |_{W_{\alpha2}}=-{\rm id} |_{W_{\alpha2}}$.

This definition makes sense, because the switch of the order of two 
subgroups $W_{\alpha1},W_{\alpha2}$ results in the change of sign 
of the involution $\varphi_\alpha$.}
\label{twodescrE}
\end{prop}

{\em Proof.} It is true for the standard divisor, and the rest follows
from the transitivity of the action of the group $G$. \hfill$\Box$

We can describe divisors $F_i$ in the same fashion, because there is 
also an involution in ${\rm Sp(4,\ZZ)}$,
whose fixed points on ${\cal H}$ are exactly the matrices $\pmatrix{x&y\cr
y&x\cr}$ that form the standard divisor $F_0.$

\begin{prop}
{ Divisors $F_i$ are in one-to-one $G$-equivariant correspondence
with conjugates in ${\rm Sp(4,\ZZ/n\ZZ)/\{\pm 1\}}$ of the involution
$$\pm\pmatrix{0&1&0&0\cr 1&0&0&0\cr 0&0&0&1\cr 0&0&1&0\cr}.$$}
\label{indexinvF}
\end{prop}

{\em Proof.} It is completely analogous to the proof of \ref{indexinvE}.
\hfill$\Box$

Now we are going to discuss the geometry of the divisor $F$.

\begin{prop}
{ Divisors $F_i$ are smooth surfaces of general type if $n$ is 
sufficiently big. Moreover, ${\rm dim}H^0(F_i,K_{F_i}) > 0$.}
\label{gentypeF}
\end{prop}

{\em Proof.} Because $F_\alpha$ is an irreducible component
of the set of fixed points of an involution on $X$, it is a smooth surface.
The finite part of $F_\alpha$ is isomorphic to
the quotient of ${\cal H}^1/\Gamma_1(2n)\times{\cal H}^1/\Gamma_1(2n)$
by the diagonal action of the group $\Gamma_1(n)/\Gamma_1(2n)$, where
${\cal H}^1$ is the usual upper half plane, and $\Gamma_1(n)$ is the principal
congruence subgroup of ${\rm Sl(2,\ZZ)}.$ This can be shown by direct calculation,
using an element of ${\rm Sp(4,\RR)}$ that maps a matrix 
$\pmatrix{x&y\cr y&x\cr}$ to the matrix $\pmatrix{x-y&0\cr 0&x+y\cr}$.
As a result, $F_\alpha$ admits a finite morphism to
$({\cal H}^1/\Gamma_1(n))^2$, which is of general type and has global 
$2$-forms, if $n$ is sufficiently big. \hfill$\Box$

The divisor $F+D$ does not have simple normal crossings.

\begin{prop}
{ There are exactly $n$ divisors $F_\gamma$ on $X_n$ that contain any
given curve $l_{\alpha\beta}=D_\alpha\cap D_\beta$.}
\label{nFl}
\end{prop}

{\em Proof.} We assume that the line $l_{\alpha\beta}$ is standard. Let us
consider the involution that fixes all points of the divisor $F_i$.
It fixes all points of the line $l_{\alpha\beta}$. This implies that 
the matrix of this involution equals
$$\pm\pmatrix{0&1&0&b\cr 1&0&-b&0\cr 0&0&0&1\cr 0&0&1&0\cr}.$$
We can lift these involutions to $\Gamma(1)$ so that they map
$\pmatrix{x&y\cr y&z\cr}$ to $\pmatrix{z+b&y\cr y&x-b\cr}$. The corresponding
divisors $F_i$ are $\pmatrix{x&y\cr y&x-b\cr}.$ The number of 
$\Gamma(n)$-inequivalent divisors of this form is equal to $n$.
\hfill$\Box$

\begin{prop}
{ If a divisor $F_\gamma$ contains a line $l_{\alpha\beta}$,
then $c_1(F_\gamma)l_{\alpha\beta}=-2$.}
\label{Fintl}
\end{prop}

{\em Proof.} The line $l_{\alpha\beta}$ may be assumed to be standard.
Calculations in the local coordinates show that the normal bindle to
$l_{\alpha\beta}$ inside $X_n$ is isomorphic to ${\cal O}(2)\oplus{\cal O}(2)$,
and the normal bundle to $l_{\alpha\beta}$ inside $F_\gamma$ is the subbundle
of the form $(x,e^{2\pi i b/n}x)$.\hfill{$\Box$}

\section{Upper bounds on the indices of subgroups of 
${\rm Sp(4,\ZZ/p^t\ZZ)}$}

This is the key section of the paper. Its purpose is to estimate
the index of the subgroup $H\subseteq {\rm Sp(4,\ZZ/n\ZZ)}$ if $H$
contains sufficiently many elements of a special type. We additionally
assume that $n=p^t$ for some prime number $p$ and integer $t$.
We fix $H$ and assume that $H\ni\pm{\bf 1}$ throughout the rest
of the section. We use the notation $[x]_p$ with $x\in {\bf R}_{\geq 1}$
for the maximum number of form $p^t,t\in {\bf N}$ that does not exceed $x$.

We first discuss subgroups that contain many elements that fix $D_i$ 
pointwise.

\begin{dfn} 
{ For any primitive vector $v$ we consider
the subgroup ${\rm Ram}_G(v)$ of $G={\rm Sp(4,\ZZ/n\ZZ)}$ that consists of 
transvections, which are operators of the form
$$r_{v,\alpha}: w \to w+\alpha\langle w,v\rangle v,~\alpha\in{\ZZ/n\ZZ}.$$ 
Because $v$ is primitive, ${\rm Ram}_G(v)\simeq {\ZZ/n\ZZ}$. We denote
$${\rm Ram}_H(v)=H\cap {\rm Ram}_G(v),~{\rm ram}_H(v)=|{\rm Ram}_H(v)|/n.$$
Clearly, ${\rm ram}_H(-v)={\rm ram}_H(v)$.}
\label{dfnramD}
\end{dfn}

\begin{rem}
{ We shall see later in Proposition \ref{ramdiv} that ${\rm Ram}_G(v_\alpha)$
is exactly the group that fixes
all points of the divisor $D_\alpha$.}
\end{rem}

\begin{prop}
{ If 
$\sum_{\pm v}{\rm ram}_H(v)\geq\epsilon \cdot \sharp (\pm v),$
then $|G:H| < 2^5\epsilon^{-2}[2^{72}\epsilon^{-42}]_p.$}
\label{boundD}
\end{prop}

{\em Proof.} We can forget about $\pm$ signs in the above proposition.
For any set $I$ of primitive vectors we define the {\it ramification 
mean}\/ of $I$ to be equal to $$(\sum_{v\in I}{\rm ram}_H(v))/|I|.$$ Clearly, 
the ramification mean never exceeds $1$. 

Among the subgroups of $V$ that are isomorphic to $({\ZZ/n\ZZ})^3$,
we choose a subgroup $V_3$, such that the ramification mean of the set 
of primitive vectors that lie in it is maximum among all such subgroups.
Any two primitive vectors are contained in the same number of subgroups
that are isomorphic to $({\ZZ/n\ZZ})^3$, so the sum of the ramification means 
among these subgroups is at least $\epsilon$ times the number of subgroups.
Hence, the ramification mean of $V_3$ is at least $\epsilon$. Analogously, 
we can choose the subgroup $V_2$ that has the maximum ramification mean 
among the subgroups with the properties

(1) $V_2\simeq ({\ZZ/n\ZZ})^2$,

(2) $V_2\subseteq V_3$,

(3) $\langle,\rangle |_{V_2}=0$.

Any two primitive vectors in $V_3$ are conjugates with respect to the
stabilizer of $V_3$ and therefore are contained in the same number of subgroups
$V_2$ that satisfy the above three properties. As a result, the ramification
mean of the set of the primitive vectors that lie in $V_2$ is also 
at least $\epsilon$. The total number of primitive vectors $v$ in $V_2$ is 
$n^2(1-p^{-2}).$ One can easily show that
at least $(\epsilon/2)n^2(1-p^{-2})$ of them have  ${\rm ram}_H(v)$ bigger than
$\epsilon/2$, because otherwise the ramification mean of $V_2$ would 
be less than $\epsilon.$ We call these vectors {\it good.}

We  may additionally assume without  loss of generality that 
$V_2={}^t(*,*,0,0)$ and $V_3={}^t(*,*,*,0).$
If $v={}^t(x,y,0,0)$, then $r_{v,1}$ has the matrix $\pmatrix{{\bf 1}&B\cr
{\bf 0}&{\bf 1}\cr}$, where $B=\pmatrix{x^2&xy\cr xy&y^2\cr}.$
Denote the group that consists of matrices $\pmatrix{{\bf 1}&B\cr
{\bf 0}&{\bf 1}\cr}$ by $G_{V_2}.$ We can prove the following
statement.

\begin{lem}
{ $|G_{V_2}:(G_{V_2}\cap H)|\leq[\epsilon^{-9}2^{14}]_p.$}
\label{V_2}
\end{lem}

{\em Proof of the lemma.}
We assume that ${\rm ram}_H({}^t(1,0,0,0))\geq \epsilon/2.$ We can do it,
because there is a primitive vector in $V_2$ with this property and
we may transform it to ${}^t(1,0,0,0)$ by an element of $G$ that stabilizes
$V_2.$ This transformation may not stabilize $V_3$, so we can not use 
this assumption in the proof of Proposition \ref{boundD}. 

At least $\epsilon n^2(1-p^{-2})/4$ good vectors ${}^t(x,y,0,0)$
satisfy ${\rm g.c.d.}(y,n)\leq [4/(\epsilon(1-p^{-2})]_p.$ Really, the number of
vectors in $V_2$ that do not satisfy this condition is at most $\epsilon 
n^2(1-p^{-2})/4.$ We pick one such vector and call it ${}^t(x_1,y_1,0,0).$

Consider the set of vectors $v={}^t(x,y,0,0)$ that 
have the following properties

(1) $v$ is good,

(2) ${\rm g.c.d.}(y,n)\leq [4/(\epsilon(1-p^{-2}))]_p$,

(3) ${\rm g.c.d.}(x_1y-y_1x,n)\leq [4/(\epsilon(1-p^{-2}))]_p$.

There are at least $\epsilon n^2(1-p^{-2})/4$ vectors that satisfy
the first two conditions and there are less than $\epsilon n^2(1-p^{-2})/4$
vectors that do not satisfy the third one. As a result, such a vector exists,
and we denote it by $v={}^t(x_2,y_2,0,0)$. 

So $H$ contains three elements of $G_{V_2}$ with the matrices
$$B=\pmatrix{\alpha^2&0\cr 0&0\cr},\pmatrix{\alpha^2x_1^2&\alpha^2x_1y_1\cr 
\alpha^2x_1y_1&\alpha^2y_1^2\cr},\pmatrix{\alpha^2x_2^2&\alpha^2x_2y_2\cr
\alpha^2x_2y_2&\alpha^2y_2^2\cr},$$
where ${\rm g.c.d.}(\alpha,n)\leq [2/\epsilon]_p.$ They generate a subgroup of $G_{V_2}$
of index equal to the greatest common divisor of $n$ and the determinant
of the corresponding three by three matrix. This is equal to 
$${\rm g.c.d.}(n,\alpha^6y_1y_2(x_1y_2-x_2y_1))\leq 
[(2/\epsilon)]_p^6[4/(\epsilon(1-p^{-2}))]^3
\leq [\epsilon^{-9}2^{14}]_p.$$
This proves the lemma. \hfill{$\Box$}

We now recall that the ramification mean of the set of vectors 
$v=~{}^t(x,y,z,0)$ is at least $\epsilon$. It implies that there are at least 
$\epsilon n^3(1-p^{-3})/2$ of them that have ${\rm ram}_H(v)\geq \epsilon/2.$
There are at least $\epsilon n^3(1-p^{-3})/4$ of them that additionally
satisfy ${\rm g.c.d.}(z,n)\leq 4/(\epsilon(1-p^{-3})).$
We abuse the notations and also call such vectors good. The operator 
$r_{v,\alpha}$ that 
corresponds to a vector $v\in V_3$ and a number $\alpha$ has the matrix
$$\pmatrix{1+\alpha xz&0&-\alpha x^2&-\alpha xy\cr \alpha yz&1&-
\alpha xy&-\alpha y^2\cr \alpha z^2&0&1-\alpha xz&-\alpha yz\cr 0&0&0&1\cr}.$$

All the elements we have described so far lie inside the group 
$$G_{V_3}=\pmatrix{a&0&b&m_3\cr m_1&1&m_2&m_4\cr
 c&0&d&m_5\cr 0&0&0&1\cr}({\rm mod}n),$$$$~ad-bc=1({\rm mod}n),
 ~bm_1+m_3=am_2({\rm mod}n),~dm_1+m_5=cm_2({\rm mod}n)\}.$$
This group has a natural projection $\lambda$ to the ${\rm Sl(2,\ZZ/n\ZZ)}$
defined by the entries $a,b,c,d.$ Our next step is to show that
the images of elements of $H$ generate a subgroup of ${\rm Sl(2,\ZZ/n\ZZ)}$
of bounded index.

We have at our disposal the matrices $\pmatrix{1+\alpha xz&-\alpha x^2\cr 
\alpha z^2&1-\alpha xz\cr}$, as well as $\pmatrix{1&\beta\cr 0&1\cr}$
with ${\rm g.c.d.}(\beta,n)\leq[\epsilon^{-9}2^{14}]_p.$ 
Here we use the estimate of $\beta$ that comes from the
result of lemma \ref{V_2}. 

We fix $\alpha_0$ that satisfies
${\rm g.c.d.}(\alpha_0,n)=[2/\epsilon]_p.$ Notice that if $(1+\alpha_0 x_1z_1,
\alpha_0 z_1^2)\neq (1+\alpha_0 x_2z_2,\alpha_0 z_2^2)$, then the matrices 
$$\pmatrix{1+\alpha_0 x_1z_1&-\alpha_0 x_1^2\cr \alpha_0 z_1^2& 
1-\alpha_0x_1z_1\cr},\pmatrix{1+\alpha_0 x_2z_2&-\alpha_0 x_2^2\cr 
\alpha_0 z_2^2&1-\alpha_0 x_2z_2\cr}$$
lie in different cosets of the subgroup $\pmatrix{1&*\cr 0&1\cr}.$
Therefore, we can estimate the order of the subgroup generated by the
elements that lie in $H$ simply by multiplying 
$n/[\epsilon^{-9}2^{14}]_p$ by the number of different pairs
$(1+\alpha_0 xz, \alpha_0 z^2)$ that we are guaranteed to have.

We have at least $\epsilon n^2(1-p^{-3})/4$ pairs $(x,z)$ that correspond 
to at least one good vector ${}^t(x,y,z,0)$ and thus give rise to an 
element in $H$ of the above form. We now need to estimate the number of 
pairs $(x,z)$ that can give the same $(1+\alpha_0 xz, \alpha_0 z^2).$
The number of different $z$ that have the same $\alpha_0 z^2$
is at most $4\cdot {\rm g.c.d.}(\alpha_0 z^2,n)$, which does not exceed
$4\cdot [2/\epsilon]_p[4/(\epsilon(1-p^{-3}))]_p^2$. Once we know $z$,
the number of $x$ that give the same $1+\alpha_0 xz$ is at most 
${\rm g.c.d.}(\alpha_0 z,n)$, which is at most
$[2/\epsilon]_p[4/\epsilon(1-p^{-3})]_p.$
So the total number of pairs $(1+\alpha_0 xz, \alpha_0 z^2)$ is at least
$$(\epsilon n^2(1-p^{-3})/4)/(4[2/\epsilon]_p^2[4/(\epsilon(1-p^{-3}))]_p^3)
\geq \epsilon n^2/(2^5[2^8\epsilon^{-5}]_p).$$

This implies that the images of elements that lie in $H$ generate a 
subgroup of ${\rm Sl(2,\ZZ/n\ZZ)}$ of index at most 
$${{n^3(1-p^{-2})(1-p^{-1})} \over {(\epsilon n^2/(2^5[2^8\epsilon^{-5}]_p))\cdot
(n/[\epsilon^{-9}2^{14}]_p)}}~\leq~[2^{22}\epsilon^{-14}]_p\epsilon^{-1}2^5.$$

On the other hand, let us estimate the index of $H\cap {\rm Ker}(\lambda)$ in
${\rm Ker}(\lambda)$. We use the formula
$$
\pmatrix{1+\alpha xz&0&-\alpha x^2&-\alpha xy\cr \alpha yz&1&-\alpha 
xy&-\alpha y^2\cr \alpha z^2&0&1-\alpha xz
&-\alpha yz\cr 0&0&0&1\cr}\cdot
\pmatrix{1&0&0&b\cr 0&1&b&0\cr 0&0&1&0\cr 0&0&0&1\cr}\cdot
$$$$
\pmatrix{1+\alpha xz&0&-\alpha x^2&-\alpha xy\cr \alpha yz&1&-\alpha 
xy&-\alpha y^2\cr  \alpha z^2&0&1-\alpha xz&-\alpha yz\cr 0&0&0&1\cr}^{
\hspace{-2pt}-1}
\hspace{-4pt}\cdot
\pmatrix{1&0&0&-b-b\alpha xz\cr 0&1&-b-b\alpha xz&b\alpha z(-2y+
bz+b\alpha xz^2)\cr 0&0&1&0\cr 0&0&0&1\cr}
$$
$$
=\pmatrix{1&0&0&0\cr -b\alpha z^2&1&0&0\cr 0&0&1&b\alpha z^2\cr 0&0&0&1}
$$
to generate the subgroup of $\pmatrix{1&0&0&0\cr *&1&0&0\cr 0&0&1&*\cr 0&0&0&1}$
of index at most ${\rm g.c.d.}(\beta\alpha z^2,n),$ which we can estimate.
$${\rm g.c.d.}(\beta\alpha z^2,n)\leq [\epsilon^{-9}2^{14}]_p[2\epsilon^{-1}]_p
[4\epsilon^{-1}(1-p^{-3})]^2_p\leq [\epsilon^{-12}2^{20}]_p.$$

Because the kernel of $\lambda$ is a semidirect product of the above group
and a subgroup of $G_{V_2}$, we have
$$|{\rm Ker}\lambda:({\rm Ker}\lambda\cap H)|\leq [\epsilon^{-12}2^{20}]_p
[\epsilon^{-9}2^{14}]_p\leq [\epsilon^{-21}2^{34}]_p$$
and
$$|G_{V_3}:(G_{V_3}\cap H)|\leq [\epsilon^{-21}2^{34}]_p[\epsilon^{-14}
2^{22}]_p\epsilon^{-1}2^5 \leq [\epsilon^{-35}2^{56}]_p\epsilon^{-1}2^5.$$

There is only one more step necessary to prove this proposition. Because the
ramification mean of $V$ is at least $\epsilon$, there are at least
$\epsilon n^4(1-p^{-4})/4$ primitive vectors $v={}^t(x,y,z,t)$ that satisfy

(1) ${\rm Ram}_H(v)\geq \epsilon/2$,

(2) ${\rm g.c.d.}(t,n)\leq [4\epsilon^{-1}(1-p^{-4})^{-1}]_p.$

We continue to abuse the notations and call these vectors good.

We use the number $\alpha_0$ defined earlier and consider elements
$r_{v,\alpha_0}$ for all good vectors. They all lie in $H$, and the claim
is that they lie in $\sim n^4$ different cosets of $G:G_{V_3}.$
Indeed, all elements of the group $G_{V_3}$ fix ${}^t(0,1,0,0)$, and 
$r_{v,\alpha_0}$ pushes ${}^t(0,1,0,0)$ to ${}^t(x\alpha_0 t, 1+y\alpha_0 t,
z\alpha_0 t,\alpha_0 t^2).$ So if 
$${}^t(x\alpha_0 t, 1+y\alpha_0 t,
z\alpha_0 t,\alpha_0 t^2)\neq {}^t(x_1\alpha_0 t_1, 1+y_1\alpha_0 t_1,
z_1\alpha_0 t_1,\alpha_0 t_1^2),$$
then $r_{v,\alpha_0}$ and $r_{v_1,\alpha_0}$ lie in different cosets.

We can estimate the number of vectors that can give the same fourtuple
as follows. If we know $\alpha_0 t^2$, it leaves us with at most
$$4\cdot {\rm g.c.d.}(\alpha_0 t^2,n)\leq 4[2\epsilon^{-1}]_p
[4\epsilon^{-1}(1-p^{-4})^{-1}]_p^2\leq 4[2^6\epsilon^{-3}]_p$$ 
options for $t.$ Once we know $t$, we have at most
$({\rm g.c.d.}(\alpha_0 t,n))^3$ choices for $(x,y,z).$
This gives us a total of at most
$$ 4[2^6\epsilon^{-3}]_p \cdot ([2\epsilon^{-1}]_p
[4\epsilon^{-1}(1-p^{-4})^{-1}]_p)^3\leq 4[2^{16}\epsilon^{-7}]_p$$
different good vectors $v$ that give $r_{v,\alpha_0}$ from the same
coset. Therefore, we can estimate the number of different cosets
that have representatives in $H$ by 
$$(\epsilon n^4(1-p^{-4})/4)/(4[2^{16}\epsilon^{-7}]_p)
\geq \epsilon n^4/(2^5[2^{16}\epsilon^{-7}]_p).$$

Hence, we can estimate the order of $H$ by multiplying the
estimate on the order of its intersection with $G_{V_3}$ by
the number of cosets that it has representatives in, which gives
$$|H|\geq {n^6(1-p^{-2})(1-p^{-1})\over 2^5\epsilon^{-1}[2^{56}\epsilon^{-35}]
_p}\cdot{\epsilon n^4\over (2^5[2^{16}\epsilon^{-7}]_p)}\geq
n^{10}\cdot{\epsilon^2 2^{-5}(1-p^{-2})(1-p^{-1})\over 
[2^{72}\epsilon^{-42}]_p}.$$

Therefore, 
$$|G:H|\leq n^{10}(1-p^{-4})(1-p^{-3})(1-p^{-2})(1-p^{-1}):
(n^{10}\cdot{\epsilon^2 2^{-5}(1-p^{-2})(1-p^{-1})\over 
[2^{72}\epsilon^{-42}]_p})
$$$$
< 2^5\epsilon^{-2}[2^{72}\epsilon^{-42}]_p.
$$
\hfill{$\Box$}

\begin{rem}
{ The estimate of Proposition \ref{boundD} is probably far from
optimum.}
\label{remboundD}
\end{rem}

Now we consider subgroups that contain many elements that fix $E_i$ 
pointwise.

\begin{dfn}
{ Let $(W_{\alpha1},W_{\alpha2})$ be a pair of complementary isotropic
subgroups that corresponds to the divisor $E_\alpha$, as described
in \ref{indexE}, and $\varphi_\alpha$ be the corresponding involution
described in \ref{twodescrE}. We define ${\rm ram}_H(E_\alpha)$ to equal $1$ 
if $H\ni \varphi_\alpha$, and to equal $0$ otherwise. This definition makes 
sense because $H\ni\pm{\bf 1}$.}
\end{dfn}

\begin{rem}
{ We have shown already that $\varphi_\alpha$ fixes all points of 
$E_\alpha$.}
\end{rem}

\begin{prop}
{ If 
$\sum_\alpha {\rm ram}_H(E_\alpha) \geq \epsilon\sharp(\alpha),$
then $|G:H| < 2^7\epsilon^{-2}[2^{246}\epsilon^{-130}]_p$.}
\label{boundE}
\end{prop}

{\em Proof.} For every set of indices $I$ we define the ramification mean
of $I$ to be $\sum_{\alpha\in I}{\rm ram}_H(E_\alpha)/|I|$.

For every primitive vector $v$ we consider the set $I_v$ of indices
$\alpha$ such that $v$ is an eigenvector of $\varphi_\alpha$.
Each index $\alpha$ belongs to the same number of sets $I_v$, therefore
$$\sum_v {\rm ramif.mean}(I_v)\geq \epsilon \sharp(v).$$

Hence there are at least $(\epsilon/2)\cdot\sharp(v)$ vectors $v$ such that
the ramification mean of $I_v$ is at least $\epsilon/2$. So now we
try to estimate ${\rm ram}_H(v)$ for a vector $v$ with this property, and then
use \ref{boundD}.

We assume that $v={}^t(0,1,0,0).$

\begin{lem}
{ Involutions $\varphi_\alpha, \alpha \in I_v$ have matrices of the
form
$$\pmatrix{1&0&0&-2x\cr -2z&-1&2x&0\cr 0&0&1&-2z\cr 0&0&0&-1}.$$
The sign is chosen to satisfy $\varphi_\alpha v=-v$.}
\end{lem}

{\em Proof of the lemma.} Any involution of this kind is defined uniquely
by the choice of $W_{\alpha_2}$. Because of $\langle 
W_{\alpha1},W_{\alpha2}\rangle =0$,
the form $\langle,\rangle$ is unimodular on $W_{\alpha2}$. It implies that 
there is a basis of $W_{\alpha2}$ that consists of $v$ and
${}^t(x,0,z,1).$ The rest is just a calculation. \hfill$\Box$

We denote the involution with the matrix 
$$\pmatrix{1&0&0&-2x\cr -2z&-1&2x&0\cr 0&0&1&-2z\cr 0&0&0&-1}$$
by $\varphi_{x,z}$. We may assume without loss of generality that 
$\varphi_{0,0}\in H$. There are at least $\epsilon n^2/2$ pairs
$(x,z)$ such that $\varphi_{x,z}\in H$. We call these pairs good.
There are at least $\epsilon n^2/4$ good pairs that satisfy
${\rm g.c.d.}(z,n) \leq [4/\epsilon]_p$. We choose one of them and denote it by
$(x_1,z_1)$. There is at least one good pair $(x,z)$ such that
${\rm g.c.d.}(xz_1-zx_1,n)\leq[2/\epsilon]\cdot {\rm g.c.d.}(z,n)$.
Then $ {\rm g.c.d.}(xz_1-zx_1,n)\leq [8\epsilon^{-2}]_p.$ We denote
this pair by $(x_2,z_2).$

It is a matter of calculation to check that 
$$(\varphi_{x_1,z_1}\varphi_{0,0}\varphi_{x_2,z_2})^2
= \pmatrix{1&0&0&0\cr 0&1&0&8(x_1z_2-x_2z_1)\cr 
0&0&1&0\cr 0&0&0&1\cr}.$$

This element lies in $H$, therefore ${\rm ram}_H(v)\geq 1/[8\epsilon^{-2}]_p.$

Because we can prove the same estimate for every vector $v$ for
which the ramification mean of $I_v$ is at least $\epsilon/2$, we get
$$\sum_v{{\rm ram}_H(v)}\geq (8\epsilon^{-2})_p^{-1}\cdot(\epsilon/2).$$

Now we use Proposition \ref{boundD} to get
$$|G:H| < 2^7\epsilon^{-2}[2^{246}\epsilon^{-130}]_p.$$
\hfill{$\Box$}

Now let us consider subgroups that contain many elements that fix lines
$D_i\cap D_j$ pointwise.

\begin{dfn}
{ Every line $l_{\alpha\beta}=D_\alpha\cap D_\beta$ is a conjugate
of the {\it standard}\/ line $l_0$, which is the intersection of the 
divisors that correspond to the $\pm$vectors $\pm{}^t(1,0,0,0)$, 
$\pm{}^t(0,1,0,0)$. We define ${\rm Ram}_G(l_0)$ to consist of matrices
$$\pmatrix{1&0&*&0\cr 0&1&0&*\cr 0&0&1&0\cr 0&0&0&1\cr} ({\rm mod}n).$$
We then define ${\rm Ram}_G(g\cdot l_0)=g\cdot {\rm Ram}_G(l_0)\cdot g^{-1}.$

It can be defined invariantly as a subgroup of all matrices that fix
both $v_\alpha$ and $v_\beta$, and also fix a pair of the isotropic
subgroups $W_1\ni v_\alpha,W_2\ni v_\beta$ that correspond to a divisor
$E_i$ that intersects $l_{\alpha\beta}$. It does not matter which $E_i$
we consider.}
\end{dfn}

\begin{rem}
{ We shall see later in Proposition \ref{jumpDD}
that ${\rm Ram}_G(l_{\alpha\beta})$ consists of
transformations that fix all points of the line $l_{\alpha\beta}$ and
do not switch the divisors $D_\alpha$ and $D_\beta$.}
\end{rem}

\begin{dfn}
{ We define ${\rm Ram}_H(l_{\alpha\beta})=H\cap {\rm Ram}_G(l_{\alpha\beta}).$ 
We define ${\rm ram}_H(l_{\alpha\beta})$ to be the maximum order of the element
of ${\rm Ram}_H(l_{\alpha\beta})$ divided by $n.$}
\end{dfn}

\begin{prop}
{ If $\sum_{\alpha\beta}{\rm ram}_H(l_{\alpha\beta})\geq \epsilon\sharp
(\alpha\beta)$, then 
$$|G:H| \leq 2^{11}\epsilon^{-2}[2^{1020}\epsilon^{-350}]_p.$$} 
\label{boundDD} 
\end{prop}

{\em Proof.} We will eventually use Proposition \ref{boundD}.
We need another definition.

\begin{lem}
{ Let $l_{\alpha\beta}$ be the line of the intersection of the divisors
$D_\alpha$ and $D_\beta.$
Then ${\rm Ram}_G(l_{\alpha\beta})={\rm Ram}_G(v_\alpha)\oplus {\rm Ram}_G(v_\beta).$}
\label{ramDD}
\end{lem}

{\em Proof of the lemma.} It is enough to consider the standard line,
for which the statement follows from the explicit matrix representations of
the three groups in question. \hfill$\Box$

\begin{dfn}
{ We define 
$${\rm ram}_H(l_{\alpha\beta}\subset D_\alpha)=
{|{\rm Ram}_H(l_{\alpha\beta})|\over |{\rm Ram}_H(l_{\alpha\beta})\cap {\rm Ram}_G(v_\alpha)|
\cdot n }$$
If $\alpha,\beta$ are standard, then ${\rm ram}_H(l_{\alpha\beta}\subset 
D_\alpha)$ is the inverse of the minimum  ${\rm g.c.d.}(a,n)$ for 
$$\pmatrix{1&0&a&0\cr 0&1&0&c\cr 0&0&1&0\cr 0&0&0&1\cr}\in H.$$}
\end{dfn}

We notice that 
$${\rm ram}_H(l_{\alpha\beta})\leq {\rm max}\{{\rm ram}_H(l_{\alpha\beta}\subset 
D_\alpha),{\rm ram}_H(l_{\alpha\beta}\subset D_\beta)\}.$$ 

The usual argument shows that at least $(\epsilon/6)\cdot\sharp|D_\alpha|$
of divisors $D_\alpha$ obey the following property

(1) at least $(\epsilon/6)\cdot\sharp (l_{\alpha\beta}\subset 
D_\alpha)$ of lines $l_{\alpha\beta}$ that are contained in it
have ${\rm ram}_H(l_{\alpha\beta}\subset D_\alpha)\geq \epsilon/6.$

We call these divisors good. Now our goal is to prove that
every good divisor $D_\alpha$ has sufficiently big ${\rm ram}_H(v_\alpha).$

We may assume without loss of generality that $D_\alpha=D_0$ is standard.
We may also assume that the arrangement of lines in $D_0$ over the standard point on
the Satake compactification contains at least $(\epsilon/6)\cdot n$
of the lines with ${\rm ram}_H(l_{0\beta}\subset D_0)\geq 
\epsilon/6.$ Divisors $D_\beta$ that intersect $D_0$ over
the standard point of the Satake compactification correspond to
$\pm$vectors of the form $\pm{}^t(1,b,0,0)$, see \ref{indexDD}.
It implies, that there are at least $(\epsilon/6)\cdot n$ numbers $b$
such that 
$$ H\ni\pmatrix{1&0&a_0&ba_0\cr 0&1&ba_0&*\cr 0&0&1&0\cr
 0&0&0&1\cr}$$
where ${\rm g.c.d.}(a_0,n) \leq [(\epsilon/6)^{-1}]_p$, and $*$ is an unknown
number.

We can choose $b_1$ and $b_2$ that give us the above elements in $H$
and additionally satisfy ${\rm g.c.d.}(b_1-b_2,n)\leq [6/\epsilon]_p.$
Then we can divide one such element by another to get 
$$H\ni\pmatrix{1&0&0&(b_1-b_2)a_0\cr
0&1&(b_1-b_2)a_0&*\cr 0&0&1&0\cr 0&0&0&1\cr}=\pmatrix{1&0&0&x\cr
0&1&x&*\cr 0&0&1&0\cr 0&0&0&1\cr}.$$
We can estimate ${\rm g.c.d.}(x,n)\leq [6^2\epsilon^{-2}]_p.$
We denote the above element by $\rho.$

Now we wander away from the standard point on the Satake 
compactification. All other divisors $D_\beta$ that intersect $D_0$
correspond to the $\pm$vectors $\pm{}^t(d,e,f,0)$ with $(d,f)\neq(0,0)(p)$.
This also follows from Proposition \ref{indexDD}. At least $(\epsilon/6) \cdot
n^3(1-p^{-2})$ of lines $l_{0\beta}$ satisfy ${\rm ram}_H(l_{0\beta}\subset 
D_0)\geq \epsilon/6.$ Therefore, at least one of them satisfies
additionally ${\rm g.c.d.}(f,n)\leq [6\epsilon^{-1}(1-p^{-2})^{-1}]_p.$
It implies, that $H$ contains an element $\rho_1$ of the form

$$\pmatrix{1+
 dfa_0&0&-d^2a_0&-dea_0\cr efa_0&1&-dea_0&-ea_0^2+c\cr 
f^2a_0&0&1-dfa_0&-efa_0\cr 0&0&0&1\cr}$$
with ${\rm g.c.d.}(f,n)\leq [6\epsilon^{-1}(1-p^{-2})^{-1}]_p.$

One can calculate that
$$\rho_1\rho\rho_1^{-1}\rho^{-1}\rho_1\rho^{-1}\rho_1^{-1}\rho=
\pmatrix{1&0&0&0\cr 0&1&0&-2x^2f^2a\cr 0&0&1&0\cr 0&0&0&1},$$
which implies 
$${\rm ram}_H(v_0)\geq [2\cdot(6^2\epsilon^{-2})^2
\cdot(6\epsilon^{-1}(1-p^{-2})^{-1})^2\cdot(\epsilon/6)^{-1}]_p^{-1}\geq
[2^{19}\epsilon^{-7}]_p^{-1}.$$

As a result,
$\sum_{\pm v}{\rm ram}_H(v)\geq [2^{19}\epsilon^{-7}]_p^{-1} 
(\epsilon/6) \cdot \sharp (\pm v).$
By \ref{boundD}, it implies
$$|G:H| \leq 2^{11}\epsilon^{-2}
[2^{1020}\epsilon^{-350}]_p.$$
\hfill$\Box$

We also need to deal with subgroups that contain many elements that fix 
$F_i$ pointwise.

\begin{dfn}
{ Let $\psi_\alpha$ be the involution that corresponds to the 
divisor $F_\alpha$ as described in \ref{indexinvF}.
We define ${\rm ram}_H(F_\alpha)$ to equal $1$ if $H\ni \psi_\alpha$,
and to equal $0$ otherwise.}
\end{dfn}

\begin{rem}
{ We have shown already that $\psi_\alpha$ fixes all points of 
$F_\alpha$.}
\end{rem}

\begin{prop}
{ If 
$\sum_\alpha {\rm ram}_H(F_\alpha) \geq \epsilon\sharp(\alpha),$
then $|G:H| \leq 2^{13}\epsilon^{-2}[2^{1722}\epsilon^{-702}]_p$.}
\label{boundF}
\end{prop}

{\em Proof.} There are at least $(\epsilon/2)\sharp(\alpha\beta)$ lines
$l_{\alpha\beta}$ such that at least $\epsilon n/2$ of divisors
$F_\gamma$ that contain $l_{\alpha\beta}$ are ramification divisors.
We call these lines good. Our goal is to estimate ${\rm ram}_H(l_{\alpha\beta})$
for a good line $l_{\alpha\beta}$.

We may assume that $l_{\alpha\beta}$ is the standard line.
If it is good, then the group $H$ contains at least $\epsilon n/2$ elements
of the form 
$$\varphi_b=\pmatrix{0&1&0&b\cr 1&0&-b&0\cr 0&0&0&1\cr 0&0&1&0\cr}.$$
There are two elements $\varphi_{b_1}$ and $\varphi_{b_2}$ in $H$ such that
${\rm g.c.d.}(n,b_1-b_2)\leq [2\epsilon^{-1}]_p$. The matrix of the element
$\varphi_{b_1}\varphi_{b_2}$ is equal to
$$\pmatrix{1&0&b_1-b_2&0\cr 0&1&0&b_2-b_1\cr 0&0&1&0\cr 0&0&0&1}.$$
Therefore, ${\rm ram}_H(l_{\alpha\beta})\geq [2\epsilon^{-1}]_p^{-1}$.
As a result,
$$\sum_{\alpha\beta}{\rm ram}_H(l_{\alpha\beta})\geq 
\epsilon2^{-1}[2\epsilon^{-1}]_p^{-1}\sharp(\alpha\beta).$$
Proposition \ref{boundDD} gives $|G:H| \leq 2^{13}\epsilon^{-2}
[2^{1722}\epsilon^{-702}]_p$. \hfill$\Box$

Finally, we will get an index estimate for subgroups such that their quotient 
varieties have bad singularities at the images of $D_i\cap D_j\cap D_k$.
This is the most delicate calculation of the whole paper. We need some
preliminary definitions.

\begin{dfn}
{ Let $P_{\alpha\beta\gamma}$ be the point of the intersection of 
three infinity divisors $D_\alpha,D_\beta$, and $D_\gamma$.
Define $${\rm Ram}_G(P_{\alpha\beta\gamma})={\rm Ram}_G(v_\alpha)\oplus
{\rm Ram}_G(v_\beta)\oplus {\rm Ram}_G(v_\gamma).$$
If $P$ is the standard point, that is the one that corresponds to
$v_\alpha={}^t(1,-1,0,0),v_\beta={}^t(1,0,0,0),v_\gamma={}^t(0,1,0,0),$
then this group consists of matrices
$$\pmatrix{1&0&a&b\cr 0&1&b&c\cr 0&0&1&0\cr 0&0&0&1\cr}({\rm mod}n).$$
As usual, we define 
${\rm Ram}_H(P_{\alpha\beta\gamma})=H\cap {\rm Ram}_G(P_{\alpha\beta\gamma})$.}
\label{ramDDD}
\end{dfn}

\begin{dfn}
{ Consider the singularity at the image of $P_{\alpha\beta\gamma}$
in the quotient of a neighborhood of $P_{\alpha\beta\gamma}$ 
by the group ${\rm Ram}_H(P_{\alpha\beta\gamma})$.
We define ${\rm mult}_H(P_{\alpha\beta\gamma})$ to be the multiplicity
of this singular point.}
\label{multDDD}
\end{dfn}

\begin{prop}
{ If $\,\sum\nolimits^* {\rm mult}_H{P_i}\geq \epsilon \sharp(i)$,
where $\sum\nolimits^*$ means taking one point $P_{\alpha\beta\gamma}$ per
orbit
of the action of the group $H$, then
$|G:H|\leq 2^{69}\epsilon^{-34}[2^{11170}\epsilon^{-5950}]_p$.}
\label{boundDDD}
\end{prop}

{\em Proof.} For each point $P_{\alpha\beta\gamma}$ we define
$\delta(H,P_{\alpha\beta\gamma})$ as a number $\delta$ defined in \ref{appdelta}
for the group ${\rm Ram}_H(P_{\alpha\beta\gamma})$ acting in the tangent space at
$P_{\alpha\beta\gamma}$. Notice that there is a natural choice of 
coordinates $(x_1,x_2,x_3)$ in a neighbourhood of $P_{\alpha\beta\gamma}$, such that the
weights of an element $h\in {\rm Ram}_H(P_{\alpha\beta\gamma})$ are determined
using ${\rm Ram}_G(P_{\alpha\beta\gamma})={\rm Ram}_G(v_\alpha)\oplus
{\rm Ram}_G(v_\beta)\oplus {\rm Ram}_G(v_\gamma)$. 
Then $\delta(H)$ is defined as $(1/n) {\rm min}_{l\neq 0}(l_1+l_2+l_3)$, where minimum 
is taken over all $H$-invariant monomials
$x_1^{l_1}x_2^{l_2}x_3^{l_3}$.

First of all, we rewrite the condition of the proposition in terms of
$\delta(H,P_{\alpha\beta\gamma})$. By \ref{appmult},
${\rm mult}_H P_{\alpha\beta\gamma} \leq n^3\delta(H,P_{\alpha\beta\gamma})
/|{\rm Ram}_H(P_{\alpha\beta\gamma})|$. Therefore,
$$\sum_{P_{\alpha\beta\gamma}}\delta(H,P_{\alpha\beta\gamma})\geq
\sum_{P_{\alpha\beta\gamma}}n^{-3}|{\rm Ram}_H(P_{\alpha\beta\gamma})|{\rm mult}_H
P_{\alpha\beta\gamma}$$
$$\geq \sum\nolimits^*(6n^3)^{-1}|H|{\rm mult}_H P_{\alpha\beta\gamma} \geq 
(6n^3)^{-1}\epsilon|H|\cdot|G:H|=\epsilon\sharp(P_{\alpha\beta\gamma}).$$

For every isotropic subgroup $V_2\simeq({\ZZ/n\ZZ})^2$ in $V$ we consider
the set of the points $P_{\alpha\beta\gamma}$ with $v_\alpha,v_\beta,
v_\gamma\in V_2$. Geometrically, these are the points that lie over certain
cusp points of the Satake compactification, see \ref{indexptsSatake}.
There are at least $(\epsilon/2)\sharp(V_2)$ of these subgroups that have
$$\sum_{v_\alpha,v_\beta,v_\gamma\in V_2}\delta(H,P_{\alpha\beta\gamma})
\geq (\epsilon/2)\sharp(v_\alpha,v_\beta,v_\gamma\in V_2).$$
We call these subgroups good. We are going to prove that if $V_2$ is a good
isotropic subgroup, then 
$$\sum_{v_\alpha,v_\beta\in V_2}{\rm ram}_H(l_{\alpha\beta})\geq
\epsilon_1(\epsilon)\sharp(v_\alpha,v_\beta\in V_2),$$
and then use Proposition \ref{boundDD}.

We assume without loss of generality that $V_2={}^t(*,*,0,0)$, and
$\delta(H,P_0)\geq (\epsilon/2)$, where $P_0$ is the standard point.
Notice that ${\rm Ram}_G(P_{\alpha\beta\gamma})$ and ${\rm Ram}_H(P_{\alpha\beta\gamma})$
do not depend on the point $P_{\alpha\beta\gamma}$, provided
$v_\alpha,v_\beta,v_\gamma\in V_2$. We denote these groups by $G_1$ and
$H_1$ respectively. The group $G_1$ is described in Definition \ref{ramDDD}.
We are dealing with points $P_{\alpha\beta\gamma}$ obtained from
the standard one by the action of elements of type $\pmatrix{A&0\cr 0&A\cr}$,
where $A\in {\rm Gl(2,\ZZ/n\ZZ)}$. Although the group $H_1$ is the same for all
$P_{\alpha\beta\gamma}$, its action in the tangent spaces depends on 
$P_{\alpha\beta\gamma}$. It is the same as the action in the tangent space
to the standard point $P_0$ of the group $AH{}^tA,~A\in{\rm Gl(2,\ZZ/n\ZZ)}$,
if we think of $G_1$ as the group of symmetric $2\times 2$ matrices.

We define $\epsilon_1$ by the formula
$$
\sum_{v_\alpha,v_\beta\in V_2}{\rm ram}_H(l_{\alpha\beta})=
\epsilon_1\sharp(v_\alpha,v_\beta\in V_2).$$

There is a line $l_{\alpha\beta}$ such that 
${\rm ram}_Hl_{\alpha\beta}\leq \epsilon_1$. It implies that the group 
$H_2=(H_1+[\epsilon_1^{-1}]_pG_1)/[\epsilon_1^{-1}]_pG_1$ is cyclic.
When we pass from $H_1$ to $H_1+[\epsilon_1^{-1}]_pG_1$, the numbers
$\delta$ do not decrease. Hence, 
$$\sum_{v_\alpha,v_\beta,v_\gamma\in V_2}\delta(H_1+[\epsilon_1^{-1}]_pG_1,
P_{\alpha\beta\gamma})\geq (\epsilon/2)\sharp(v_\alpha,v_\beta,v_\gamma\in 
V_2).$$ 
This is equivalent to 
$$\sum_{A\in {\rm Gl(2,\ZZ/n\ZZ)}}\delta(AH_1{}^tA+[\epsilon_1^{-1}]_pG_1,P_0)
\geq (\epsilon/2)\sharp(A).$$

Let $H_2$ be generated by $B=\pmatrix{a&b\cr b&c\cr}$. One can show that
$\delta(AH{}^tA+[\epsilon_1^{-1}]_pG_1,P_0)$ equals $\delta(({\ZZ/n\ZZ})
{\bar A}{\bar B}{}^t{\bar A},{\bar P}_0)$, where $n$ is replaced by 
$[\epsilon_1^{-1}]_p$ and bars means reduction 
${\rm mod}[\epsilon_1^{-1}]_p$. Therefore, 
$$\sum_{C\in {\rm Gl(2,\ZZ/}[\epsilon_1^{-1}]_p\ZZ)}
\delta((\ZZ/[\epsilon_1^{-1}]_p\ZZ)CB{}^tC,{\bar P_0})\geq
(\epsilon/2)\sharp(C).$$

Because of the result of Proposition \ref{appfinmany},
there are at most $2^{10}\epsilon^{-8}[2^{12}\epsilon^{-5}]_p$ different
matrices $CB{}^tC$ up to 
proportionality that give $\delta((\ZZ/[\epsilon_1^{-1}]_p\ZZ)
CB{}^tC,{\bar P_0})\geq (\epsilon/4)\sharp(C)$. This implies that
the orbit $CB{}^tC({\rm mod~proportionality})$ of the action of the group 
${\rm Gl(2,\ZZ/}[\epsilon_1^{-1}]_p\ZZ)$ has length at most 
$2^{12}\epsilon^{-9}[2^{12}\epsilon^{-5}]_p$.

However, we can estimate this length by looking at matrices $C=
\pmatrix{t&0\cr 0&1\cr}$. They give $CB{}^tC=\pmatrix{t^2a&tb\cr tb&c\cr}$,
and so length of the orbit is at least $[\epsilon_1^{-1}]_p(1-p^{-1})
/{\rm g.c.d.}(bc,[\epsilon_1^{-1}]_p).$ Because we have assumed that $\delta
(\ZZ/[\epsilon_1^{-1}]_p\ZZ)B,{\bar P_0})\geq (\epsilon/2)$,
we have ${\rm g.c.d.}(b,n)\leq[2\epsilon^{-1}]_p$ and ${\rm g.c.d.}(c,n)\leq
[4\epsilon^{-1}]_p$. Really, the weights of $B$ are $(-b,a+b,c+b)$,
and if ${\rm g.c.d.}[c,n]$ is greater than $[4\epsilon^{-1}]_p$, then 
we get $\delta\geq (\epsilon/2)$ because of the invariant monomial of the form
$(x_1x_3)^{[\epsilon_1^{-1}]_p/[4\epsilon^{-1}]_p}$.

As a result, the length of the orbit is at least $[\epsilon_1^{-1}]_p
(1-p^{-1})/[8\epsilon^{-2}]_p$, and we have
$[\epsilon_1^{-1}]_p(1-p^{-1})/[8\epsilon^{-2}]_p\leq 2^{12}\epsilon^{-9}[2^{12}
\epsilon^{-5}]_p$ and $\epsilon_1\geq 2^{-28}\epsilon^{16}$.

We now may use the result of Proposition \ref{boundDD} with
$(2^{-28}\epsilon^{16})(\epsilon/2)$ in place of $\epsilon$.
Thus, $|G:H|\leq 2^{69}\epsilon^{-34}[2^{11170}\epsilon^{-5950}]_p$.
\hfill$\Box$

\section{Singularities of ${\cal H}/H$}

It is easy to describe all elements of finite order in $\Gamma(2)$
by means of the following proposition.

\begin{prop}
 { Any nonidentity element of finite order in $\Gamma(2)/\{\pm 1\}$
 is conjugate in $\Gamma(1)/\{\pm 1\}$ to the element with the matrix
 $$\varphi_0=\pmatrix{1&0&0&0\cr0&-1&0&0\cr0&0&1&0\cr0&0&0&-1\cr}.$$}
\label{einv}
\end{prop}

{\em Proof.} Denote the matrix of this element by $\varphi=\pmatrix
 {A&B\cr C&D\cr}.$ Because $\Gamma(4)$ is torsion free and
 $\varphi^2\in\Gamma(4)$, we obtain $\varphi^2=1$. Hence the following
 equalities hold $$A~{}^tB=B~{}^tA,~~C~{}^tD=D~{}^tC,~~A~{}^tD-B~{}^tC=1$$
 $$A=~{}^tD,~~B=-{}^tB,~~C=-{}^tC.$$Really, the first three equalities hold
 for all symplectic matrices, and they imply  $\varphi^{-1}=\pmatrix
 {{}^tD&-{}^tB\cr -{}^tC&{}^tA\cr}$, so $\varphi^{-1}=\varphi$ gives the
 last three ones. Six equalities together show that
 $$\varphi=\pmatrix{a_1&a_2&0&b\cr a_3&a_4&-b&0\cr0&c&a_1&a_3\cr
 -c&0&a_2&a_4\cr}$$ with $(a_1+a_4)b=(a_1+a_4)c=(a_1+a_4)a_2=
 (a_1+a_4)a_3=0,~a_1^2+a_2a_3-bc=a_4^2+a_2a_3-bc=1.$ Hence if
 $\varphi\neq1$, then $a_1+a_4=0$, so
 $$\varphi=\pmatrix{a_1&a_2&0&b\cr
 a_3&-a_1&-b&0\cr0&c&a_1&a_3\cr -c&0&a_2&-a_1\cr}$$
 with $a_1^2+a_2a_3
 -bc=1,~(a_1-1),a_2,a_3,b,c\equiv 0{\rm mod}(2)$.

We need to prove that any matrix with these properties is conjugate to
 $\varphi_0$. The vector spaces ${\rm Ker}(\varphi-1)$
 and ${\rm Ker}(\varphi+1)$ are orthogonal, so we should simply find
 four integer vectors $e_1,...,e_4$ that obey $\varphi(e_i)=(-1)^{i+1}
 e_i$ and $\langle e_2,e_4\rangle =\langle e_1,e_3\rangle =1$. 
 Because of symmetry, it is enough
 to find $e_1$ and $e_3$. Let us denote $d={\rm g.c.d.}(b/2,a_3/2,(a_1-1)/2).$
 There holds $\alpha b/2+\beta(a_1-1)/2+\gamma(-a_3/2)=d$ for some
 integers $\alpha,\beta,\gamma.$ Now we simply put 
$$e_1=\pmatrix{b/2d\cr 0\cr a_3/2d\cr(1-a_1)/2d\cr},~ 
e_3=\alpha\pmatrix{0\cr-b/2\cr(a_1+1)/2\cr a_2/2\cr}+\beta
 \pmatrix{a_2/2\cr(1-a_1)/2\cr c/2\cr 0\cr}+\gamma\pmatrix{(a_1+1)
 /2\cr a_3/2\cr0\cr-c/2\cr}$$
 and check the required conditions by direct calculation. \hfill$\Box$

\begin{dfn}
 { Let $H$ be a subgroup of finite index in $\Gamma(1)$. We call
 $E_i$ or $F_j$ a {\em ramification divisor} iff $H$ contains the
 involution that fixes all points of the divisor. Because of 
 the results of \ref{indexinvE} and \ref{indexinvF}, $E_\alpha$ is a
 ramification divisor iff ${\rm ram}_H(E_\alpha)=1$, and similarly
for $F_\beta$.}
\end{dfn}

We are interested in singularities of ${\cal H}/H$. They occur
 at the images of the points of ${\cal H}$ that have nontrivial
 stabilizers in $H$. The goal of the rest of this section is to prove
 the following statement.

\begin{prop}
{ Singularities of the images of the points $\xi\in{\cal H}$ that
 do not lie in ramification divisors $E_i$ or $F_j$ are
 canonical. Points that do lie in ramification divisors have
 solvable stabilizers of order at most $72$. We refer to 
 \ref{defcan} for the definitions of canonical and terminal 
singularities. }
\label{finsing}
\end{prop}

{\em Proof. }There are two possibilities: $\xi\in\cup E_i$ and
 $\xi\notin\cup E_i$.

{\em Case 1. } $\xi\notin\cup E_i$. The stabilizer of $\xi$ in $\Gamma(2)$
 equals $\{\pm 1\}$ because of Proposition
 \ref{einv} and the definition of $E_i$. We
 consider the quotient of ${\cal H}$ by the action of $\Gamma(2)$.
 It is the smooth part of the singular quartic $V$ defined by the equation
 $(\sum x_i^2)^2=4\sum x_i^4$ in coordinates $(x_1:...:x_6,~\sum x_i=0)$
 of ${\bf P}^4$, see \cite{Geer}. The group $\Gamma(1)/\Gamma(2)\simeq
 \Sigma_6$ acts on $V$ by the permutations of the coordinates $x_i.$ The
 stabilizer  $\xi$ in $\Gamma(1)$ equals that of the image of $\xi$ in $V$ in
 the group $\Sigma_6.$ Moreover, locally their actions are the same, so
 the resulting quotient singularities are isomorphic. Therefore, we need to 
 study fixed points of $\Sigma_6$-action on $V$.

\begin{lem}
{ A point  $\xi\notin\cup E_i$ with a nontrivial stabilizer in $\Gamma(1)$
 either lies in $\cup F_j$ or has the image in $V$ of type
 $\sigma(0:\theta:\theta^2:\theta^3:\theta^4:1),~\theta=\exp(2\pi i/5),
 ~\sigma\in\Sigma_6.$}
\label{fixedsigma}
\end{lem}

{\em Proof of the lemma.} Denote by $x=(x_1:...:x_6)$ the image
 of $\xi$ in $V$. We may assume that the stabilizer of $x$ contains one
 of the permutations
 $$(1,2);(1,2)(3,4);(1,2)(3,4)(5,6);(1,2,3);(1,2,3)(4,5,6);(1,2,3,4,5).$$
 Let us calculate the sets of fixed points of these permutations that
lie in $V$.

Case (1,2). We have $(x_1,x_2,...,x_6)=\lambda(x_2,x_1,...,x_6).$
 If $\lambda=-1$, then $x=(-1:1:0:0:0:0)$, but this point does not lie in $V$.
 Hence $\lambda=1$. The set defined by "$x_1=x_2$" constitutes an irreducible divisor on
 $V$, so it is the closure of the image of some submanifold of dimension two
 in ${\cal H}$.

Case (1,2)(3,4). We have $(x_1,x_2,x_3,x_4,x_5,x_6)=\lambda(x_2,x_1,x_4,
 x_3,x_5,x_6).$ If $\lambda=-1$, then $x_1=-x_2,~x_3=-x_4,~x_5=x_6=0$.
 The equality $(\sum x_i^2)^2=4\sum x_i^4$ implies that
 $x_1=x_3$ or $x_1=x_4$, so $x\in {\rm Sing}(V)$, see \cite{Geer}. If $\lambda=1$,
 then $x$ lies in the divisor "$x_1=x_2$".

Case (1,2),(3,4),(5,6). We have $(x_1,x_2,x_3,x_4,x_5,x_6)=\lambda
 (x_2,x_1,x_4,x_3,x_6,x_5).$ If $\lambda=-1$, then
 $x_1=-x_2,~x_3=-x_4,~x_5=-x_6.$ Equality $(\sum x_i^2)^2=4\sum x_i^4$
 leads to $(x_1+x_3+x_5)\cdot(x_1+x_3-x_5)\cdot(x_1-x_3+x_5)\cdot(
 -x_1+x_3+x_5)=0.$ Each of these linear equations implies that $x$
 lies in the image of $\cup E_i$, see \cite{Geer}.
 If $\lambda=1$, then $x_1=x_2,~x_3=x_4,~x_5=x_6$ so $x\in {\rm Sing}(V)$.

Case (1,2,3). We have $(x_1,x_2,x_3,x_4,x_5,x_6)=\lambda(x_2,x_3,x_1
 ,x_4,x_5,x_6).$ If $\lambda\neq1$, then $x_1+x_2+x_3=0$, so $x$ lies in
 the image of $\cup E_i$. If $\lambda=1$, then $x_1=x_2=x_3$, and $x$ lies
 in the divisor "$x_1=x_2$".

Case (1,2,3)(4,5,6). We have $(x_1,x_2,x_3,x_4,x_5,x_6)=\lambda(x_2,
 x_3,x_1,x_5,x_6,x_4).$ If $\lambda=1$, then
 $x=(1:1:1:-1:-1:-1)\notin V$. Otherwise, $x_1+x_2+x_3=0$, and $\xi\in\cup E_i$.

Case (1,2,3,4,5). We have $(x_1,x_2,x_3,x_4,x_5,x_6)=\lambda(x_2,x_3
 ,x_4,x_5,x_1,x_6).$ If $\lambda=1$, then
 $x=(1:1:1:1:1:-5)\notin V$. Otherwise $x=\sigma(0:\theta:
 \theta^2:\theta^3:\theta^4:1),~\theta=\exp(2\pi i/5),~\sigma\in\Sigma_6.$

The above calculation shows that there is only one up to $\Sigma_6$-action
 divisor on $V$ with a nontrivial stabilizer of a generic closed point.
 On the other hand, the images of $F_j$ on $V$ obey this condition.
 Therefore the images of $F_j$ are the conjugates of the divisor "$x_1=x_2$",
 which proves the lemma. \hfill$\Box$

\begin{rem}
{ As a corollary of this lemma, codimension one components
of the ramification locus of the map from ${\cal H}/\Gamma_n$
to ${\cal H}/H$ can only be divisors $E_\alpha$ and $F_\beta$. 
Moreover, ramification occurs
iff $E_\alpha(F_\beta)$ is a ramification divisor as defined
above, and in this case the only nontrivial element that preserves all points
of the divisor is the corresponding involution. Of course, when we consider
the Igusa compactifications, we may have ramification at infinity divisors.}
\label{ramfin}
\end{rem}

Let us come back to the proof of \ref{finsing}. We try to estimate 
the singularity
 at the image of the point $\xi\notin \cup E_i$ under the quotient map ${\cal
 H}\to{\cal H}/H.$ The group $\Gamma(2)/\{\pm 1\}$ acts freely on
 ${\cal H}-\cup E_i$, so we can work in terms of the image point
 $x\in V-{\rm Sing}(V)$ and the group ${\rm Stab}^H\xi\cdot\Gamma(2)/\Gamma(2)$,
 because these quotient singularities are isomorphic.
 There exists a useful criterion that enables one to find out
 whether the quotient singularity is canonical, see \cite{Reid}.
 In particular, it is always canonical, if the image of the group in ${\rm
 Gl}(T_x)$ lies in ${\rm Sl}(T_x)$. We use these facts
 extensively.

 First of all we consider the case $x=\sigma(0:\theta:...:1)$. Then
 either ${\rm Stab}^H\xi\cdot\Gamma(2)/\Gamma(2)=1$ or ${\rm Stab}^H\xi\cdot
 \Gamma(2)/\Gamma(2)={\ZZ/5\ZZ}.$ A direct calculation of the weights of
 the generator in the tangent space and the criterion of \cite{Reid}
 show that the quotient singularity is terminal, hence canonical.

Now let us consider other points $x=(x_1:x_2:x_3:x_4:x_5:x_6)\in V-{\rm Sing}(V)$.
 The group $S={\rm Stab}^H\xi\cdot\Gamma(2)/\Gamma(2)$ contains no transpositions,
 because $\xi$ does not belong to any ramification divisors $F_i$.
 The proof of \ref{fixedsigma} shows that $S$ does not contain
 permutations of types $(*,*)(*,*)(*,*),(*,*,*)(*,*,*)$, and
 $(*,*,*,*,*)$. As a result, $S$ consists of permutations of types $(*,*)
 (*,*),(*,*,*),(*,*,*,*)$, and $(*,*)(*,*,*,*)$ only.
 Calculations similar to those of \ref{fixedsigma} show that if the group $S$
 contains $(*,*)(*,*,*,*)$, then $\xi\in \cup E_i$. Moreover, if it
 contains a permutation of type $(*,*,*,*)$ and the proportionality
 coefficient $\lambda$ does not equal $1$, then  $\xi\in \cup E_i$.
 Notice (see the proof of \ref{fixedsigma}) that the proportionality
 coefficients of elements of the group $S$ of types $(*,*)(*,*)$ and $(*,*,*)$
 must also equal $1$. All these restrictions on the group $S$ imply that it
 consists of even permutations, and all proportionality coefficients are
 equal to $1$. Therefore, the group $S$ acts in the tangent space
 of $x$ by matrices from ${\rm  Sl}$. The criterion of M. Reid
 shows that the quotient singularity is canonical.

{\em Case 2.} $\xi\in \cup E_i$. all divisors $E_i$ are conjugates,
so we may assume that $\xi$ is represented by a diagonal matrix. Different 
$E_i$ do not intersect, so the stabilizer $S$ of $\xi$ in $\Gamma(1)$
 is a subgroup of $${\rm Stab}(\Delta)=\{\pmatrix{
 a&0&b&0\cr0&a_1&0&b_1\cr c&0&d&0\cr0&c_1&0&d_1\cr}\cup\pmatrix{
 0&a_1&0&b_1\cr a&0&b&0\cr 0&c_1&0&d_1\cr c&0&d&0\cr},~ad-bc=a_1
 b_1-c_1d_1=1\}.$$
 Point $\xi$ may be transformed by the group ${\rm Stab}(\Delta)$ to the point
 $\pmatrix{x_0&0\cr0&z_0\cr}$ with $|{\rm Re}(x_0)|\leq 1/2,~|x_0|\geq
 1,~|{\rm Re}(z_0)|\leq 1/2,~|z_0|\geq 1$. Without any loss of generality
 one may consider points of this type only. The stabilizer of
 the general such point in $\Gamma(1)/\{\pm 1\}$ equals ${\ZZ/2\ZZ}.$
 It is generated by the involution of Proposition \ref{einv}.
 If this element is in $H$, then $\Delta$ is a ramification divisor
 by our definition. The order of the stabilizer can increase in the
 following curves and points (we have used the symmetry between $x$
and $z$)
 $$\pmatrix{x&0\cr0&x\cr},\pmatrix{i&0\cr0&x\cr},\pmatrix{\rho&0\cr0&x\cr}
 ,\pmatrix{i&0\cr0&\rho\cr},\pmatrix{i&0\cr0&i\cr},
 \pmatrix{\rho&0\cr0&\rho\cr}.$$

Let us check all these cases.

Case $\pmatrix{x&0\cr0&x\cr}$. Because $\Delta$ is not a ramification
 divisor, the order of ${\rm Stab}^H\xi$ is at most two, so the quotient
 singularity is canonical.

Case $\pmatrix{i&0\cr0&x\cr}$. We get $| {\rm Stab}^H\xi|=1$ by the same
 argument.

Case $\pmatrix{\rho&0\cr0&x\cr}$. In this case either $| {\rm Stab}^H\xi|=1$
 or $ {\rm Stab}^H\xi$ is generated by the element of order $3$ whose action
 in the tangent space of $\xi$ has determinant $1$.

Case $\pmatrix{i&0\cr0&\rho\cr}$. The argument is the same as in the
 previous case.

Case $\pmatrix{i&0\cr0&i\cr}$. The stabilizer of $\xi$ in $\Gamma(1)
 /\{1,-1\}$ is generated by the images of elements of $\Gamma(1)$
 with matrices
 $$\varphi=\pmatrix{0&1&0&0\cr 1&0&0&0\cr 0&0&0&1\cr 0&0&1&0\cr},
 ~\alpha=\pmatrix{1&0&0&0\cr 0&0&0&1\cr 0&0&1&0\cr 0&-1&0&0\cr},
 ~\beta=\pmatrix{0&0&1&0\cr 0&1&0&0\cr -1&0&0&0\cr 0&0&0&1\cr}.$$
 The relations are $\alpha\beta=\beta\alpha,~\alpha^2=\beta^2,~ \varphi\alpha
 =\beta\varphi,~\varphi^2=\alpha^4=\beta^4=1.$ The order of the group
 is $16$.

The point $\xi$ does not lie in the ramification divisors, so ${\rm Stab}^H\xi$
 does not contain  any conjugates of $\varphi$. As in the above cases,
 ${\rm Stab}^H\xi$ does not contain $\alpha^2$ either. We can also employ the
 following simple statement: if $s^2=1$ for all $s\in {\rm Stab}^H\xi$, then
the quotient singularity is canonical. In our case it implies that if the 
quotient singularity is not canonical, then the group ${\rm Stab}^H\xi$ contains an
 element of order $4$. All these conditions on  ${\rm Stab}^H\xi$ together
 hold iff this group is generated by a conjugate of $\varphi\beta.$
 A direct calculation of the weights in the tangent space completes the
 argument.

Case $\pmatrix{\rho&0\cr0&\rho\cr}.$ In this case it is possible to check
 that  the condition "$\xi$ does not belong to any ramification divisors" 
implies that ${\rm Stab}^H\xi$ acts in the tangent space of
 $\xi$ by matrices with determinant $1$.

To finish the proof of Theorem \ref{finsing}, we only need to check that
 stabilizers of all points of ${\cal H}$ are solvable groups whose
orders are at most $72$. It can be done using the description of 
${\cal H}/\Gamma(2)$ as the smooth part of the singular quartic. I skip
the details, because this number is clearly bounded and only slightly
affects the constant in the final result.
\hfill$\Box$

\section{Finiteness Theorem for subgroups $H\supseteq \Gamma(p^t)$}

We assume that $n=p^t$ throughout this section. We denote the subgroup
$H\supseteq \Gamma(n)$ and the quotient $H/\Gamma(n)$ by the same letter, which
should not lead to a confusion. The Igusa compactifications of ${\cal H}
/\Gamma(n)$ and ${\cal H}/H$ are denoted by $X$ and $Y$. The quotient map
$X\to X/H=Y$ is denoted by $\mu$.

We start by pulling the problem from $Y$ to $X$.

\begin{dfn}
{ Let $\pi:Z\to Y$ be a desingularization of $Y$. Denote by $-1+\delta$
the minimum discrepancy of $Y$, see \ref{discrep}.  Because of 
\ref{qulog}, $\delta$ is a positive rational number.}
\end{dfn}

\begin{dfn}
{ Let $m$ be a sufficiently divisible number, so that $mK_Y$ is
a Cartier divisor on $Y$. The vector space $H^0(Y,mK_Y-{\rm mlt})$ is
defined as the space of global sections $s$ of the coherent subsheaf of
${\cal O}_Y(mK_Y)$ that consists of sections that lie in
$m_y^{m(1-\delta)}({\cal O}_Y(mK_Y))_y$ for all noncanonical singular points
$y\in Y$.}
\end{dfn}

\begin{rem}
{ We assume $m$ to be sufficiently divisible whenever it
is necessary. We also omit ${\cal O}$ in the notations of the 
space of global sections, unless it can lead to a misunderstanding.}
\end{rem}

\begin{prop}
{ ${\rm dim}H^0(Z,mK_Z)\geq {\rm dim}H^0(Y,mK_Y-{\rm mlt}).$}
\label{ZtoY}
\end{prop}

{\em Proof.} The pullbacks $\pi^*s$ vanish with the multiplicity at least
$m(1-\delta)$ along exceptional divisors with negative discrepancies.
Hence we can define an injective linear map from $H^0(Y,mK_Y-{\rm mlt})$
to $H^0(Z,mK_Z)$. \hfill$\Box$

\begin{dfn}
{ Denote by $H^0(Y,mK_Y-{\rm mlt}^0)$ the space of global
sections $s\in H^0(Y,mK_Y)$ that satisfy $s\in m_y^{m(1-\delta)}({\cal 
O}_Y(mK_Y))_y$ for all noncanonical singular points of $Y$ except for
the images of points $P_{\alpha\beta\gamma}$ that are triple intersections
of infinity divisors on $X$.}
\end{dfn}

Clearly, $H^0(Y,mK_Y-{\rm mlt}^0))\supseteq H^0(Y,mK_Y-{\rm mlt})$.

\begin{prop}
{ If $|G:H| > 2^{953}[2^{165870}]_p$, then
${\rm dim} H^0(Y,mK_Y-{\rm mlt}^0)-{\rm dim}H^0(Y,mK_Y-{\rm mlt})
\preceq_{m\to\infty} 2^{-8}3^{-6}5^{-1}m^3|G:H|.$}
\label{excludeDDD}
\end{prop}

{\em Proof.} When $m\to\infty$, the codimension we are trying
to estimate grows no faster than
$(\sum_{Q\in Y}{\rm mult}_Q)(m^3/6)$, where $\sum_{Q\in Y}{\rm mult}_Q$ is the
sum over all points $Q$ in the image of $\cup P_{\alpha\beta\gamma}$,
and ${\rm mult}_Q$ is the multiplicity of the local ring of $Y$ at $Q$.
We want to relate it to the statement of Proposition \ref{boundDDD}.
We need an easy lemma.

\begin{lem}
{ Let $P_{\alpha\beta\gamma}=D_\alpha\cap D_\beta \cap D_\gamma$
be a point on $X$, such that $\mu (P_{\alpha\beta\gamma})=Q$.
Then ${\rm mult}_Q\leq 6^3 {\rm mult}_H(P_{\alpha\beta\gamma})$ with  
${\rm mult}_H(P_{\alpha\beta\gamma})$ defined in \ref{multDDD}.}
\end{lem}

{\em Proof of the lemma.} Every element of $G$ that fixes 
$P_{\alpha\beta\gamma}$ permutes the triple of the $\pm$vectors
$(\pm v_\alpha, \pm v_\beta, \pm v_\gamma)$. Hence the subgroup in ${\rm Stab}^H(P)$
of the elements that induce trivial permutations is a normal subgroup
of order at most $6$. One can show that this subgroup coincides with 
${\rm Ram}_H(P_{\alpha\beta\gamma})$ by the explicit matrix calculation for the
standard triple $v_\alpha={}^t(0,1,0,0),v_\beta={}^t(-1,1,0,0), v_\gamma=
{}^t(1,0,0,0)$.

Therefore, the singularity of $Y$ at $Q$ can be obtained as the quotient
of the singularity of $X/{\rm Ram}_H(P_{\alpha\beta\gamma})$ by the group of
 order at most $6$. Its multiplicity can be estimated by means of
Proposition \ref{klem}. \hfill$\Box$

As a result of this lemma, the codimension we are trying to estimate grows
no faster than $m^36^2\sum\nolimits^*{\rm mult}_H(P_{\alpha\beta\gamma})$, where
one takes 
one point $P_{\alpha\beta\gamma}$ per orbit of $H$. By \ref{boundDDD} with
$\epsilon=2^{-26}$, it grows no faster than $2^{-8}
3^{-6}5^{-1}m^3|G:H|$, if $|G:H| > 2^{953}[2^{165870}]_p$.
\hfill$\Box$

\begin{prop}
{ Let $L_Y$ be a divisor of the modular form of weight $1$ on $Y$.
Then ${\rm dim}H^0(Y,mL_Y)$ grows as $2^{-7}3^{-6}5^{-1}m^3|G:H|$.}
\label{growthL}
\end{prop}

{\em Proof.} It can be derived, for instance, from the formula for
${\rm dim}H^0(X,mL_X)$ and $\oplus_mH^0(Y,mL_Y)=(\oplus_mH^0(X,mL_X))^H$.
\hfill$\Box$

\begin{prop}
{ If ${\rm dim} H^0(Y,m(K_Y-L_Y)-{\rm mlt}^0)\neq 0$
for sufficiently big $m$ and $|G:H| > 
2^{953}[2^{165870}]_p $, then the variety $Y$ is of general type.}
\label{getridofDDD}
\end{prop}

{\em Proof.} We get
$${\rm dim}H^0(Z,mK_Z)\geq {\rm dim}H^0(Y,mK_Y-{\rm mlt})$$
$$\geq {\rm dim}H^0(Y,mK_Y-{\rm mlt}^0)-2^{-8}3^{-6}5^{-1}|G:H|m^3$$
$$\succeq 
{\rm dim}H^0(Y,mL_Y)-2^{-8}3^{-6}5^{-1}|G:H|m^3\sim
2^{-8}3^{-6}5^{-1}|G:H|m^3.$$
\hfill$\Box$

We shall eventually prove that if $|G:H|$ is big, then
${\rm dim}H^0(Y,m(K_Y-L_Y)-{\rm mlt}^0)\neq 0$ for big $m$.

\begin{dfn}
{ Let $R$ be the ramification divisor of the morphism $\mu$.
We define $H^0(X,mK-mR-mL-{\rm mlt}^0)$ to be the space
of global sections of ${\cal O}_X(m(K_X-R-L_X))$ that satisfy
certain vanishing conditions. Namely, we require their germs to lie
in $m_x^{m\cdot k({\rm Stab}^H(x))} {\cal O}_X(m(K_X-R-L_X))_x$ for all
points $x\in X$ whose images in $Y$ have noncanonical singularities,
except for $x=P_{\alpha\beta\gamma}.$ Here $k({\rm Stab}^H(x))$ is defined 
according to remark \ref{k}.}
\end{dfn}

\begin{prop} 
{ If $|G:H|> 2^{953}[2^{165870}]_p$ and 
${\rm dim}H^0(X,mK-mR-mL-{\rm mlt}^0)\neq 0$ for some $m>0$,
then the variety $Y$ is of general type.}
\label{alreadyonX}
\end{prop}

{\em Proof.} Because of \ref{klem}, all $H$-invariant elements
of $H^0(X,mK-mR-mL-{\rm mlt}^0)$ can be pushed down to 
elements of $H^0(Y,m(K_Y-L_Y)-{\rm mlt}^0).$ Notice, that
$\mu^*(m(K_Y-L_Y))=m(K_X-R-L_X)$, and $m\delta$ is dropped from the
vanishing conditions to compensate for the constant $N$ from \ref{klem}.
One can multiply the $H$-conjugates of any section to get an $H$-invariant
one, so if ${\rm dim}H^0(X,mK-mR-mL-{\rm mlt}^0)\neq 0$ for some $m>0$,
then ${\rm dim}H^0(Y,m(K_Y-L_Y)-{\rm mlt}^0)\neq 0$ for big $m$,
and Proposition \ref{getridofDDD} finishes the proof.\hfill$\Box$

Let us now describe ramification divisors and points with noncanonical
images.

\begin{prop} 
{ The ramification divisor $R$ equals 
$\sum_\alpha(n\cdot {\rm ram}_H(v_\alpha)-1)D_\alpha +
\sum_\alpha {\rm ram}_H(E_\alpha)E_\alpha +\sum_\alpha {\rm ram}_H(F_\alpha)F_\alpha.$}
\label{ramdiv}
\end{prop}

{\em Proof.} We know by \ref{ramfin} that the ramification divisor in
the finite part is equal to $\sum_\alpha {\rm ram}_H(E_\alpha)E_\alpha+\sum_\alpha 
{\rm ram}_H(F_\alpha)F_\alpha$. We only need to show that the group of
elements of $G$ that fix all points of the divisor $D_\alpha$ is
exactly $\pm {\rm Ram}_G(v_\alpha)$. It can be done explicitly in coordinates
for the standard divisor $D_0.$ \hfill$\Box$

\begin{prop}
{ If $x\in D_\alpha$, but $x\notin \cup (D_\alpha\cap D_\beta)$,
then ${\rm Stab}^H(x)/(\pm {\rm Ram}_H(v_\alpha))$ is a group of  order at most $6$.}
\label{jumpD}
\end{prop}

{\em Proof.} We only need to consider the standard divisor $D_0$.
It is the universal elliptic curve with level $n$ structure.
It can be shown that the group ${\rm Stab}^G_{D_0}$ acts on it by
a combination of modular transformations of the base, additions
of points of order $n$ in the fibers, and the involution
$a\to -a$ of the fibers. The order $6$ can be reached for the point
$x$ on the curve with complex multiplication, such that $x$ satisfies
$2n\cdot x=0$, and all other stabilizers are even smaller. I skip
the details, because a different bound here would only slightly
affect the final estimate. \hfill$\Box$

\begin{prop}
{ If $x\in l_{\alpha\beta}= D_\alpha\cap D_\beta$,
but $x\notin \cup P_{\alpha\beta\gamma}$,
then the order of the group
${\rm Stab}^H(x)/\pm {\rm Ram}_H(l_{\alpha\beta})$ is at most
$4$.}
\label{jumpDD}
\end{prop}

{\em Proof.} We may assume that $l_{\alpha,\beta}=l_0$ is the standard line.
The group ${\rm Stab}^G(l_0)$ contains a subgroup of index $2$ of elements
that preserve both $\pm{}^t(1,0,0,0)$ and $\pm{}^t(0,1,0,0).$
It in turn contains a subgroup of
index $2$ that consists of matrices $\pm \pmatrix{{\bf 1}&B\cr
{\bf 0}&{\bf 1}\cr}, ~B=\pmatrix{a&b\cr b&c\cr}.$
One can show using the explicit coordinate on $l_0$, that 
if $b$ is nonzero, then the action of this element has no fixed points
on $l_0$, except for the points of triple intersection of the 
infinity divisors, which finishes the proof. \hfill$\Box$

\begin{prop}
{ If $\,|G:H|> 2^{953}[2^{165870}]_p$ and
$${\rm dim}H^0(m(K-L)-m\sum_\alpha n\cdot {\rm ram}_H(v_\alpha) 7 D_\alpha
-m\sum_\alpha {\rm ram}_H(E_\alpha) 73 E_\alpha$$
$$-m\sum_\alpha {\rm ram}_H(F_\alpha) 73 F_\alpha
 -m\sum_{\alpha\beta}n\cdot {\rm ram}_H(l_{\alpha\beta}) 4 l_{\alpha\beta})
\neq 0$$
for some $m>0$, then the variety $Y$ is of general type.}
\label{73}
\end{prop}

{\em Proof.} We know from Proposition \ref{finsing} that the
points in the finite part, that do not
lie in the ramification divisors $E_\alpha$ or $F_\beta$,
do not contribute to ${\rm mlt}^0.$ Therefore,
$$m\sum_\alpha {\rm ram}_H(E_\alpha) 73 E_\alpha+m\sum_\alpha {\rm ram}_H(F_\alpha)
73 F_\alpha \geq {\rm mlt}+mR$$
in the finite part. This inequality, strictly speaking, is the inclusion of
the sheaves of ideals. Analogously, Propositions \ref{jumpD} and
\ref{jumpDD} show that 
$$m\sum_\alpha n\cdot {\rm ram}_H(v_\alpha) 7 D_\alpha\geq mR_D+{\rm mlt}^0$$
on $D$ away from $\cup(D_\alpha\cap D_\beta)$, and 
$$m\sum_{\alpha\beta}n\cdot {\rm ram}_H(l_{\alpha\beta}) 4 l_{\alpha\beta}
\geq {\rm mlt}^0$$
on $\cup(D_\alpha\cap D_\beta)$ away from points $P_{\alpha\beta\gamma}$.
Then it remains to use Proposition \ref{alreadyonX}.
\hfill$\Box$

\begin{prop}
{ If the variety $Y$ is not of general type, then at least
one of the following inequalities holds true.

(1) $$|G:H|\leq  2^{953}[2^{165870}]_p$$

(2) $$ {\rm dim}H^0(m(K-L))-{\rm dim}H^0(m(K-L)-m\sum_\alpha n\cdot {\rm ram}_H(v_\alpha) 
7 D_\alpha)$$
$$\geq ((1/6)c_1(K-L)^3m^3)/5$$

(3) $$ {\rm dim}H^0(m(K-L))-{\rm dim}H^0(m(K-L)-m\sum_\alpha {\rm ram}_H(E_\alpha)
73 E_\alpha)$$
$$\geq ((1/6)c_1(K-L)^3m^3)/5$$

(4) $$ {\rm dim}H^0(m(K-L))-{\rm dim}H^0(m(K-L)-m\sum_\alpha {\rm ram}_H(F_\alpha)
 73 F_\alpha)$$
$$\geq ((1/6)c_1(K-L)^3m^3)/5$$

(5) $$ {\rm dim}H^0(m(K-L)-m\sum_\alpha n\cdot {\rm ram}_H(v_\alpha) 7 D_\alpha
-m\sum_\alpha {\rm ram}_H(E_\alpha) 73 E_\alpha$$
$$-m\sum_\alpha {\rm ram}_H(F_\alpha) 73 F_\alpha)
-{\rm dim}H^0(m(K-L)-m\sum_\alpha n\cdot {\rm ram}_H(v_\alpha) 7 D_\alpha$$
$$-m\sum_\alpha {\rm ram}_H(E_\alpha) 73 E_\alpha
-m\sum_\alpha {\rm ram}_H(F_\alpha) 73 F_\alpha$$
$$ -m\sum_{\alpha\beta}n\cdot {\rm ram}_H(l_{\alpha\beta}) 4 l_{\alpha\beta})
\geq ((1/6)c_1(K-L)^3m^3)/5$$}
\label{splitcases}
\end{prop}

{\em Proof.} If (2),(3), and (4) are all false, then 
$${\rm dim}H^0(m(K-L)-m\sum_\alpha n\cdot {\rm ram}_H(v_\alpha) 7 D_\alpha
-\sum_\alpha {\rm ram}_H(E_\alpha) 73 E_\alpha$$
$$-\sum_\alpha {\rm ram}_H(F_\alpha) 73 F_\alpha)\succeq (2/5)(1/6)c_1(K-L)^3m^3.$$
Really, ${\rm dim}H^0(m(K-L)$ grows like $(1/6)c_1(K-L)^3m^3$, because $K-L$ is
ample for big $n$, and $E_\alpha, F_\beta, D_\gamma$ are different divisors.
Hence, if (1) and (5) are also false, then Proposition 
\ref{73} proves that the variety $Y$ is of general type. \hfill$\Box$

Our next goal is to show that each of the statements (2)-(5)
implies that $|G:H|$ is less than some constant. 
We use results of Yamazaki \cite{Yamazaki} and statements of
Section 3. 

\begin{prop}
{ If $$ {\rm dim}H^0(m(K-L))-{\rm dim}H^0(m(K-L)-m\sum_\alpha n\cdot {\rm ram}_H(v_\alpha) 
7 D_\alpha)$$
$$\geq ((1/6)c_1(K-L)^3m^3)/5,$$
then $|G:H|< 2^{41}[2^{828}]_p$.}
\label{5.D}
\end{prop}

{\em Proof.} First of all, we get 
$$ {\rm dim}H^0(m(K-L))-{\rm dim}H^0(m(K-L)-m\sum_\alpha n\cdot {\rm ram}_H(v_\alpha) 
7 D_\alpha)$$
$$\leq \sum_\alpha({\rm dim}H^0(m(K-L))-{\rm dim}H^0(m(K-L)-7mn\cdot {\rm ram}_H(v_\alpha)
D_\alpha)).$$
The standard exact sequences associated to $D_\alpha\subset X$ allow us
to estimate that 
$${\rm dim}H^0(m(K-L))-{\rm dim}H^0(m(K-L)-7mn\cdot {\rm ram}_H(v_\alpha)
D_\alpha)$$
$$\leq \sum_{j=0}^{7mn\cdot {\rm ram}_H(v_\alpha)-1}
{\rm dim}H^0(D_\alpha,m(K-L)-jD_\alpha)$$
$$=\sum_{j=0}^{7mn\cdot {\rm ram}_H(v_\alpha)-1}
{\rm dim}H^0(D_\alpha,m(K-L)+(2j/n)(L+E)).$$

The divisor $L+E$ is nef on $X$, because $L$ is nef,
divisors $E_i$ are disjoint, and $(L+E)E_i=0$. The divisor $K-L$
is ample on $X$, if $n$ is sufficiently big. Therefore, we may use the 
Riemann-Roch
formula to calculate ${\rm dim}H^0(D_\alpha,m(K-L)+(2j/n)(L+E))$. Because we
are only interested in the coefficient of $m^3$, as $m\to \infty$,
we only need to take into account the term
$(1/2)c_1(m(K-L)+(2j/n)(L+E))^2c_1(D_\alpha)$.
When $j$ grows, this intersection number grows, therefore
$${\rm dim}H^0(m(K-L))-{\rm dim}H^0(m(K-L)-7mn\cdot {\rm ram}_H(v_\alpha)
D_\alpha)$$
$$\leq 7mn\cdot {\rm ram}_H(v_\alpha) (1/2) m^2c_1(K-L+14{\rm ram}_H(v_\alpha)(L+E))^2
c_1(D_\alpha)$$
$$\leq m^3 \sharp(v_\alpha)^{-1} {\rm ram}_H(v_\alpha) 
(7n/2)c_1(K-L+14 (L+E))^2c_1(D).$$

Hence, if the condition of the proposition is true, then 
$$\sharp(v_\alpha)^{-1}\sum_\alpha {\rm ram}_H(v_\alpha)\geq
105^{-1}c_1(K-L)^3/(c_1(K-L+14(L+E))^2c_1(nD)).$$
The right hand side can be calculated using the formulas of Yamazaki for
the intersection numbers of the divisors $D,L,K$, and $E$.
It is bigger than $2^{-18}$ if $n$ is sufficiently big, which we may assume
without loss of generality. Therefore, by the result of Proposition 
\ref{boundD}, $|G:H|< 2^{41}[2^{828}]_p$.\hfill$\Box$

\begin{prop}
{ If  $$ {\rm dim}H^0(m(K-L))-{\rm dim}H^0(m(K-L)-m\sum_\alpha {\rm ram}_H(E_\alpha)
73 E_\alpha)$$
$$\geq ((1/6)c_1(K-L)^3m^3)/5,$$
then $|G:H|< 2^{53}[2^{3236}]_p$.}
\label{5.E}
\end{prop}

{\em Proof.} Analogously to the proof of \ref{5.D}, we estimate
$$ {\rm dim}H^0(m(K-L))-{\rm dim}H^0(m(K-L)-m\sum_\alpha {\rm ram}_H(E_\alpha)
73 E_\alpha)$$
$$\leq \sum_\alpha {\rm ram}_H(E_\alpha) \sum_{j=0}^{73m-1}
{\rm dim}H^0(E_\alpha,m(K-L)-jE_\alpha)$$
$$=\sum_\alpha {\rm ram}_H(E_\alpha) 
\sum_{j=0}^{73m-1} {\rm dim}H^0(E_\alpha,m(K-L)+jL)$$
$$\preceq
\sharp(E_\alpha)^{-1}\sum_\alpha {\rm ram}_H(E_\alpha) (73/2)m^3 
c_1(K+72L)^2c_1(E).$$
Therefore, 
$$\sharp(E_\alpha)^{-1}\sum_\alpha {\rm ram}_H(E_\alpha)
\geq 73^{-1}15^{-1}c_1(K-L)^3/(c_1(K+72L)^2c_1(E))>2^{-23}.$$
Then Proposition \ref{boundE} tells us that 
$|G:H|< 2^{53}[2^{3236}]_p$.\hfill$\Box$

\begin{prop} 
{ If $$ {\rm dim}H^0(m(K-L))-{\rm dim}H^0(m(K-L)-m\sum_\alpha {\rm ram}_H(F_\alpha)
 73 F_\alpha)$$
$$\geq ((1/6)c_1(K-L)^3m^3)/5,$$
then $|G:H|< 2^{73}[2^{22782}]_p$.}
\label{5.F}
\end{prop}

{\em Proof.} As in the proof of \ref{5.E}, we estimate 
$$ {\rm dim}H^0(m(K-L))-{\rm dim}H^0(m(K-L)-m\sum_\alpha {\rm ram}_H(F_\alpha)
73 F_\alpha)$$
$$\leq \sum_\alpha {\rm ram}_H(F_\alpha) \sum_{j=0}^{73m-1}
{\rm dim}H^0(F_\alpha,m(K-L)-jF_\alpha).$$
Unfortunately, the geometry of $F$ is more complicated than that of $E$,
and we do not have a nice formula like $(L+E_\alpha)E_\alpha=0$. We can
get away with it by using the adjunction formula together with the Proposition
\ref{gentypeF}.

We can estimate 
$$\sum_\alpha {\rm ram}_H(F_\alpha) \sum_{j=0}^{73m-1}
{\rm dim}H^0(F_\alpha,m(K-L)-jF_\alpha)$$
$$\leq \sum_\alpha {\rm ram}_H(F_\alpha) \sum_{j=0}^{73m-1}
{\rm dim}H^0(F_\alpha,m(K-L)+jK-jK_{F_\alpha})$$
$$\leq \sum_\alpha {\rm ram}_H(F_\alpha) \sum_{j=0}^{73m-1}
{\rm dim}H^0(F_\alpha,m(K-L)+jK)$$
$$\leq \sharp(F_\alpha)^{-1}\sum_\alpha {\rm ram}_H(F_\alpha)
(73/2)m^3c_1(74K-L)^2c_1(F).$$
Therefore, 
$$\sharp(F_\alpha)^{-1}\sum_\alpha {\rm ram}_H(F_\alpha)\geq
73^{-1}15^{-1}c_1(K-L)^3/(c_1(74K-L)^2c_1(F)).$$
We need to have some upper bound on $c_1(74K-L)^2c_1(F).$
To do this, we recall the proof of Proposition \ref{finsing},
where we have shown that the images of the divisors $F_\alpha$ on the singular
quartic $V$ have form $x_i=x_j$. The product $\prod_{i\neq j}(x_i-x_j)^2$
is invariant under the permutations of the coordinates, so it defines a
modular form of weight $60$, that vanishes on $F$. Here we use the fact
that the coordinates of ${\bf P}^4$ are given by the modular forms of weight $2$,
see \cite{Geer}.
As a result, $c_1(74K-L)^2c_1(F)\leq 60c_1(74K-L)^2c_1(L)$,
and we can estimate 
$\sharp(F_\alpha)^{-1}\sum_\alpha {\rm ram}_H(F_\alpha)>2^{-30}$.

Now Proposition \ref{boundF} implies that $|G:H|< 2^{73}[2^{22782}]_p$.
\hfill$\Box$

\begin{prop}
{ If 
$$ {\rm dim}H^0(m(K-L)-m\sum_\alpha n\cdot {\rm ram}_H(v_\alpha) 7 D_\alpha
-m\sum_\alpha {\rm ram}_H(E_\alpha) 73 E_\alpha$$
$$-m\sum_\alpha {\rm ram}_H(F_\alpha) 73 F_\alpha)
-{\rm dim}H^0(m(K-L)-m\sum_\alpha n\cdot {\rm ram}_H(v_\alpha) 7 D_\alpha$$
$$-m\sum_\alpha {\rm ram}_H(E_\alpha) 73 E_\alpha
-m\sum_\alpha {\rm ram}_H(F_\alpha) 73 F_\alpha$$
$$ -m\sum_{\alpha\beta}n\cdot {\rm ram}_H(l_{\alpha\beta}) 4 l_{\alpha\beta})
\geq ((1/6)c_1(K-L)^3m^3)/5,$$
then $|G:H|<2^{65}[2^{10470}]_p$.}
\label{5.DD}
\end{prop}

{\em Proof.} Denote 
$$L_1=K-L-7\sum_\alpha n\cdot {\rm ram}(v_\alpha)D_\alpha - 73\sum_\alpha
{\rm ram}_H(E_\alpha)E_\alpha-73\sum_\alpha {\rm ram}_H(F_\alpha)F_\alpha.$$
Then the left hand side of the proposition does not exceed the sum over
all $l_{\alpha\beta}$ of
$${\rm dim}H^0(mL_1)-{\rm dim}H^0(mL_1-4mn\cdot {\rm ram}_H(l_{\alpha\beta})l_{\alpha\beta}).$$

To estimate this codimension, we consider the blow-up of the variety $X$
along the line $l_{\alpha\beta}$, which we denote by $\pi:X_1\to X$.
The normal bundle to $l_{\alpha\beta}$ is isomorphic to 
${\cal O}(2)\oplus{\cal O}(2)$. This can be checked by direct calculation.
Therefore, the exceptional divisor of $\pi$ is isomorphic to ${\bf P}^1\times
{\bf P}^1$. We get
$${\rm dim}H^0(mL_1)-{\rm dim}H^0(mL_1-4mn\cdot {\rm ram}_H(l_{\alpha\beta})l_{\alpha\beta})$$
$$={\rm dim}H^0(m\pi^*L_1)-{\rm dim}H^0(m\pi^*L_1-4mn\cdot {\rm ram}_H(l_{\alpha\beta})S)$$
$$\leq\sum_{j=0}^{4mn\cdot {\rm ram}(l_{\alpha\beta})-1}{\rm dim}H^0(S,m\pi^*L_1-jS).$$

We denote the fiber and the section of $S\to l_{\alpha\beta}$ by $f$ and
$s$ respectively and get $(m\pi^*L_1-jS)S=m\cdot c_1(L_1)l_{\alpha\beta}
\cdot f +j(2f+s)$. Hence, $H^0(S,m\pi^*L_1-jS)$ grows when $j$ grows, and
we have 
$${\rm dim}H^0(mL_1)-{\rm dim}H^0(mL_1-4mn\cdot {\rm ram}_H(l_{\alpha\beta})l_{\alpha\beta})$$
$$\leq 4mn\cdot {\rm ram}_H(l_{\alpha\beta}){\rm dim}H^0(S,mc_1(L_1)l_{\alpha\beta}\cdot
2f+4mn\cdot {\rm ram}_H(l_{\alpha\beta})\cdot j )$$
$$\leq m^3{\rm ram}_H(l_{\alpha\beta})\cdot (8n\cdot {\rm ram}_H(l_{\alpha\beta})
+c_1(L_1)l_{\alpha\beta})\cdot4n\cdot {\rm ram}_H(l_{\alpha\beta})$$
$$\leq m^3{\rm ram}_H(l_{\alpha\beta})(128n^3+16n^2c_1(L_1)l_{\alpha\beta})$$
$$\preceq_{n\to\infty}~m^3{\rm ram}_H(l_{\alpha\beta})(128n^3+
16n^2\cdot (7n\cdot 2\cdot 2+73\cdot 2\cdot n)
$$
$$=m^3{\rm ram}_H(l_{\alpha\beta})\cdot 2912n^3.$$

The number of $l_{\alpha\beta}$ is equal to $2^{-3}n^7(1-p^{-4})(1-p^{-2})$,
see \cite{Yamazaki}. Therefore, if the condition of the proposition is true,
then 
$$(\sharp(l_{\alpha\beta}))^{-1}\sum_{\alpha\beta}{\rm ram}_H(l_{\alpha\beta})
\geq {{c_1(K-L)^3
}\over{
(30\cdot 2912\cdot 2^{-3}n^{10}(1-p^{-4})(1-p^{-2}))}}>2^{-27}.$$
Now the result of Proposition \ref{boundDD} gives
$|G:H|<2^{65}[2^{10470}]_p$. \hfill$\Box$

We are now ready to prove the finiteness theorem for $H\supseteq 
\Gamma(p^t)$.

\begin{prop}
{ If $|G:H|>2^{953}[2^{165870}]_p$,
then the variety $Y$ is of general type.}
\end{prop}

{\em Proof.} We simply combine the results of Propositions
\ref{5.D}, \ref{5.E}. \ref{5.F}, \ref{5.DD}, and \ref{splitcases}.
\hfill$\Box$

\begin{prop}
{ {\bf Finiteness theorem for $H\supseteq \Gamma(p^t).$}
There are only finitely many subgroups $H\subseteq {\rm Sp(4,\ZZ)}$
of finite index that contain $\Gamma(p^t)$ for some $p$ and $t$,
such that the variety ${\cal H}/H$ is not of general type.}
\label{finthmprimary}
\end{prop}

{\em Proof.} It follows from the fact that $|G:H|$ is bounded.
\hfill$\Box$

In particular, if $p$ is sufficiently big, then for any $H$,
${\rm Sp(4,\ZZ)}\supset H\supseteq \Gamma(p^t)$ the variety $Y$
is of general type.

\section{Finiteness Theorem, the general case}

Now we no longer assume that $n$ is a power of a prime number. Our goal
is to prove that the condition $n=p^t$ can be dropped from the statement
of Proposition \ref{finthmprimary}. Our proof is the direct
generalization of the argument of \cite{Thompson}.

We first estimate prime factors of $n$.

\begin{prop}
 { If $p>3$, and
$$H\cdot\Gamma(p)=\Gamma(1),~
 H\supseteq\Gamma(mp^\alpha),~{\rm g.c.d.}(m,p)=1,$$ 
then
$H\supseteq\Gamma(m)$.}
\label{bigprime}
\end{prop}

{\em Proof. } For any group $G$ we denote its image modulo
 $\Gamma(mp^\alpha)$ by $\hat G$. We have isomorphisms
$$\hat\Gamma(1)\simeq
 \hat\Gamma(m)\times\hat\Gamma(p^\alpha),~\hat\Gamma(m)\simeq{\rm
 Sp(4,\ZZ/p^\alpha \ZZ)},~\hat\Gamma(p^\alpha)\simeq{\rm Sp(4,\ZZ/m
\ZZ)}.$$

The group ${\rm PSp(4,\ZZ/p \ZZ)}$ is simple for all $p\geq 3$.
 Because of $\hat H\cdot\hat\Gamma(p)/\hat\Gamma(p)\simeq
 {\rm Sp(4,\ZZ/p \ZZ)}$, the group $\hat H$ has a section isomorphic to
 ${\rm PSp(4,\ZZ/p \ZZ)}$. Consider the following normal subgroups of
 $\hat\Gamma(1)$.
$$\hat\Gamma(1)\supset\hat\Gamma(m)\supset\hat\Gamma(mp)\supseteq\{e\}.$$
 We easily get that $\hat H\cap\hat\Gamma(m)/\hat H\cap\hat\Gamma(mp)$
 has a section isomorphic to ${\rm PSp(4,\ZZ/p\ZZ)}$, so there holds
 $$(\hat H\cap\hat\Gamma(m))\cdot\hat\Gamma(mp)=\hat\Gamma(m)$$

Now it will suffice to prove that the last equality implies $\hat H\supseteq
 \hat\Gamma(m)$. Note that $\hat\Gamma(m)\simeq {\rm Sp(4,\ZZ/p^\alpha
 \ZZ)}$ and $\hat\Gamma(mp)\simeq {\rm Ker}({\rm Sp(4,\ZZ/p^\alpha \ZZ)}
 \to{\rm Sp(4,\ZZ/p \ZZ}))$. We denote by $K_i$ the kernels of ${\rm Sp(4
 ,\ZZ/p^\alpha \ZZ)}\to{\rm Sp(4,\ZZ/p^i \ZZ)}$
 for $i=1,...,\alpha$ and prove that $\hat H\supseteq K_i$ by the decreasing
 induction on $i$.

For $i=\alpha$ there is nothing to prove. Besides we already have the last
 step of the induction. Suppose that  $\hat H\supseteq K_i,~i>1$. To prove
 that  $\hat H\supseteq K_{i-1}$ consider $h\in\hat H\cap\hat\Gamma(m)$
 such that $$h\equiv \pmatrix{
 1&0&1&0\cr0&1&0&0\cr0&0&1&0\cr0&0&0&1}({\rm mod}~p).$$
 Clearly, $h^{p^i}\in K_i$. Besides, a simple calculation shows that for
 $p\geq 5$ $$h^{p^{i-1}}\equiv\pmatrix{1&0&p^{i-1}&0\cr0&1&0&0\cr0&0&1&0
 \cr0&0&0&1}({\rm mod}~p^i).$$ When the group $\hat\Gamma(m)$ acts on
 $K_{i-1}/K_i$ by conjugation, its subgroup $\hat\Gamma(mp)$ acts as
 identity. We have already known that $(\hat H\cap\hat\Gamma(m))\cdot
 \hat\Gamma(mp)=\hat\Gamma(m)$, so it is enough to show that
 conjugates of the element $h^{p^{i-1}}$ generate the whole group
 $K_{i-1}$ modulo $K_i$. This can be done by a direct calculation in the
 abelian group $K_{i-1}/K_i$. \hfill $\Box$

\begin{prop}
{ There exists a natural number $N$ such that if ${\cal H}/H$ is not
of general type, then $$H\supseteq \Gamma(\prod_{p_i\leq N}p_i^{n_i})$$
for some natural numbers $n_i$.}
\label{boundprime}
\end{prop}

{\em Proof.} Let $n$ be the minimum number such that $H\supseteq \Gamma(n)$.
Because of the result of \ref{bigprime}, $H\cdot \Gamma(p)\neq \Gamma(1)$
for all prime factors of $p$ of $n$ bigger than $3$. If ${\cal H}/H$ is not of
general type, then ${\cal H}/(H\cdot\Gamma_p)$ is not of general type
either, see \ref{fincov}. Because of Proposition \ref{finthmprimary}, 
there are only finitely many choices for $p.$ \hfill$\Box$

We now prove the Finiteness Theorem in full generality.
 
Define for any $H\subseteq\Gamma(1)$ and any prime $p$ the
 $p$-projection of $H$ as $H_p=\cap_1^\infty H\cdot\Gamma(p^j)$. Note
 that $H_p\supseteq H$ and $H_p\supseteq\Gamma(p^j)$ for some $j$. The
 following proposition allows us to work with $p$-projections only,
 after we have got an estimate on the primes.

\begin{prop}
 { For any given set of subgroups $G_i\supseteq\Gamma(p_i^{n_i}),
 ~i=1,...,k$, there are only finitely many subgroups $H\supseteq\Gamma
 (p_1^{\alpha_1}\cdot...\cdot p_k^{\alpha_k})$ with $H_{p_i}=G_i$.}
\label{splitprime}
\end{prop}

{\em Proof.} We can simply estimate the index of $H$ if we employ
 the fact that $\Gamma(p_i)$ are pro-$p_i$-groups.\hfill $\Box$

Now we can easily prove the Finiteness Theorem.

\begin{prop}
{ {\bf Finiteness Theorem.} There are only finitely many subgroups 
$H\subseteq {\rm Sp(4,\ZZ)}$ of finite index, such that ${\cal H}/H$
is not of general type.}
\label{fintheorem}
\end{prop}

{\em Proof.} If ${\cal H}/H$ is not of general type, then 
${\cal H}/H_p$ is not of general type either. Therefore, Proposition
\ref{finthmprimary} tells us that there are only finitely many choices
for $H_p$. By \ref{boundprime}, all prime factors of $n$ are bounded,
so Proposition \ref{splitprime} finishes the proof.
\hfill $\Box$

\section{Varieties of general type and singularities}

We first recall some standard facts about varieties of general type and 
singularities.

\begin{dfn}
{ A smooth compact algebraic variety X over $\CC$ is called a variety of
general type if there exists some constant $c>0$ such that ${\rm dim}H^0(X,{\cal
O}_X(mK_X))>cm^{{\rm dim}X}$ for all sufficiently big
(equivalent condition -- divisible by some integer $d$) positive integers $m$.
Here  $K_X$ is the canonical divisor of $X$.}
\end{dfn}

\begin{rem}
{ If $X$ and $Y$ are birational smooth compact algebraic varieties, then
${\rm dim}H^0(X,{\cal O}_X(mK_X))= {\rm dim}H^0(Y,{\cal O}_Y(mK_Y))$ for $m\geq0$. }
\end{rem}

\begin{dfn}
{ A field ${\cal K}\supset\CC$ is called a field of general type if it
is a field of the rational functions
of a smooth compact algebraic variety of general type.}
\end{dfn}

\begin{dfn}
{ An algebraic variety over $\CC$ is called a variety of general type if
its field of functions is a field of general type.}
\end {dfn}

\begin{dfn}
{ A canonical divisor $K_X$ of a normal variety $X$ is a Weil divisor 
on $X$ that coincides with a canonical divisor on $X-{\rm Sing}(X)$. The variety $X$
is called $\QQ$-$Gorenstein$ if $mK_X$ is a Cartier divisor for some integer
$m$.}
\end{dfn}

\begin{rem}
{ If the variety $Y$ is normal  $\QQ$-Gorenstein but has singularities,
then the condition "${\rm dim}H^0(Y,{\cal O}_Y(mK_Y))>cm^{{\rm dim}Y}$ for $m\to+\infty$"
does not imply by itself that $Y$ is of general type. Really, if
$\pi:Z\to Y$ is some desingularization, then there holds
$$K_Z=\pi^*(K_Y)+\sum_i\alpha_iF_i,~~\alpha_i\in\QQ$$
in the sense of equivalence of $\QQ$-Cartier divisors, where $F_i$ are
exceptional divisors of morphism $\pi$ and $\alpha_i$ are some rational numbers
called  {\it discrepancies}\/. If some $\alpha_i$ are negative,
then  ${\rm dim}H^0(Z,{\cal O}_Z(mK_Z))$ may be less than ${\rm dim}H^0(Y,{\cal
O}_Y(mK_Y))$}.
\label{discrep}
\end{rem}

\begin{dfn}
{ A normal  $\QQ$-Gorenstein variety $Y$ is said to have $log-terminal$
 singularities if for some desingularization $\pi:Z\to Y$, such that the
exceptional divisor $\sum F_i$ has simple normal crossings, all 
discrepancies are greater than $(-1)$. A singular point $y\in Y$ is called 
{\it canonical}\/ (resp. {\it terminal}\/) if the discrepancies $\alpha_i$ are 
nonnegative (resp. positive) for all $i$ such that $\pi(F_i)\ni y$.
Once satisfied for some desingularization, whose exceptional locus is a divisor
with simple normal crossings, these conditions are satisfied for any
desingularization (see \cite{CKM}).}
\label{defcan}
\end{dfn}

\begin{prop}
{ If $\mu:X\to Y$ is a finite morphism of algebraic varieties and $Y$ is of
general type, then $X$ is also of general type.}
\label{fincov}
\end{prop}

{\em Proof.} We find a surjective morphism $\hat\mu:\hat X\to\hat Y$,
where $\hat X,\hat Y$ are smooth projective birational models of $X,Y$,
and then pull back multicanonical forms.
\hfill$\Box$

The following statement is well-known.

\begin{prop}
{ {\rm (see \cite{CKM})}
Let $X$ be a smooth projective algebraic variety over $\CC$ with an
action of a finite group $G$.
Then the quotient variety $Y=X/G$ has log-terminal singularities.}
\label{qulog}
\end{prop}

Now we shall prove a simple but important technical result about quotient
singularities. Let $X$ be a projective algebraic variety with an action of a
finite solvable group $H$. Let $x$ be a (closed) point of $X$, such that $hx=x$
for all $h\in H$. Suppose we have $\{e\}=H_0\subset H_1\subset ...\subset 
H_t=H$, where $H_{i-1}$ are normal subgroups of $H_i$ and $H_i/H_{i-1}$ are 
abelian groups with exponents $k_i$. Denote $k=k_1\cdot...\cdot k_t$. 
Denote the local ring of $x$ in $X$ by $(A,\mm)$. Then $(B,\nn)=(A^H,\mm^H)$ is the 
local ring of the image of $x$ under the quotient morphism.

\begin{prop}
{ In the above setup there exists a constant $N$, which depends only on $X$
and $H$ but not on $x$, such that there holds $\mm^{kl+N}\cap B\subseteq \nn^l$ 
for all $l\geq0$.}
\label{klem}
\end{prop}

{\em Proof.} We do not suppose $X$ to be smooth, so it is enough to consider
just the case of an abelian group $H$ with $kH=0$. There exists a linearized
$H$-invariant very ample invertible sheaf ${\cal L}$ on $X$. Consider the 
corresponding closed embedding $X\to {\bf P}^{N_0}$. Because
$H$ is abelian, there exists an open $H$-invariant affine neighborhood of
$x$ with the ring $R$ equal to $\CC[1,l_1/l_0,...,l_{N_0}/l_0]/I$
where
$l_i\in H^0(X,{\cal L}),~
h(l_i)=\mu_i(h)\cdot l_i,~\forall h\in H$ and $I$ is some ideal. Moreover,
we may assume
that $f_i=l_i/l_0$ vanish at $x$, because of $Hx=x$. Hence the local ring
$(A,\mm)$ is the localization of $R$ in $p=(f_1,...,f_{N_0})$.
Because $H$ is finite, one can assume that all denominators are 
$H$-invariant. Therefore, the statement of the proposition
is equivalent to $p^{kl+N}\cap R^H\subseteq (p^H)^l$.

Each element of $p$ can be represented as a polynomial in $f_i$ with zero
constant term. Therefore, each element of $p^{ kl+N}$ can be represented as a
polynomial in $f_i$ with monomials of degree no less than $kl+N$. For any
given $f\in p^{kl+N}\cap R^H$ consider such a representation with the minimum
possible number of monomials. Then if for some monomial $g$ of this
representation and some element $h\in H$ there holds $h(g)=w\cdot g,~w\neq1$, 
then the formula $f=f\cdot w/(w-1)-h(f)/(w-1)$ allows us to reduce the number
of monomials. Hence every element $f\in p^{kl+N}\cap R^H$ is a sum of 
$H$-invariant monomials of degree at least $kl+N$.

Now we only need to prove that any $H$-invariant monomial $g=f_1^{\alpha_1}
\cdot...\cdot f_{N_0}^{\alpha_{N_0}}$ of degree at least $kl+N$ is a product
of at least $l$ $H$-invariant monomials of positive degree. It is time to
choose $N$, namely $N=k\cdot N_0$. Denote by $\gamma_i$ the maximum integers
that do not exceed $\alpha_i/k$. Then $g=f_1^{k\gamma_1}\cdot...\cdot 
f_{N_0}^{k\gamma_{N_0}}\cdot g_1$ gives the required decomposition,
because $\sum\gamma_i>\sum\alpha_i/k-N_0\geq l$.\hfill $\Box$

\begin{rem}
{ Due to the result of \cite{Hochster}, the above proposition holds
for scheme points which correspond to the subvarieties that are pointwise
$H$-invariant. I wish to thank Melvin Hochster for pointing out this 
reference.}
\label{afterklem}
\end{rem}

\begin{rem}
{ In the rest of the paper $k(H)$ for a finite solvable group $H$ denotes
the least possible value of $k$ that could be obtained in the above way.}
\label{k}
\end{rem}

The rest of the section is devoted to multiplicities of certain toric 
singularities. Somewhat unnatural choice of notation is motivated by 
the notation of Section 3.

\begin{dfn}
{ Let $G_1\simeq ({\ZZ/n\ZZ})^3$ act on $\CC^3$ according to the formula
$$(\xi_1,\xi_2,\xi_3)(x_1,x_2,x_3)=(e^{2\pi i\xi_1/n}x_1,
e^{2\pi i\xi_2/n}x_2,e^{2\pi i\xi_3/n}x_3).$$
Let $H_1$ be a subgroup of $G_1$. Define $\delta(H_1)=(1/n)
{\rm min}_{l\neq 0}(l_1+l_2+l_3)$, where the minimum is taken
among all $H_1$-invariant monomials
$x_1^{l_1}x_2^{l_2}x_3^{l_3}$.}
\label{appdelta}
\end{dfn}

\begin{prop}
{ The multiplicity of the local ring of $C^3/H_1$ at zero is at most
$n^3\delta(H_1)/|H_1|.$}
\label{appmult}
\end{prop}

{\em Proof.} The exponents of the $H_1$-invariant monomials form a semigroup,
which we denote by $K$. One can show that the multiplicity is equal to
$vol(\RR^n_{>0}-conv(K-\{0\}))/|H_1|$, where the volume is normalized
to be equal one on the basic tetrahedron. This result does not seem to
be stated explicitly anywhere in the literature, but its proof is completely
analogous to the calculation of \cite{Teissier} of multiplicities
of the ideals in the polynomail ring that are generated by monomials.
On the other hand, this set is contained in the set
$$conv((l_1,l_2,l_3),(0,0,n),(0,n,0),(0,0,0))\cup...$$
$$...\cup
conv((l_1,l_2,l_3),(0,n,0),(n,0,0),(0,0,0)),$$
which has volume $n^3\delta(H_1)$. \hfill$\Box$

\begin{rem}
{ Our results on the multiplicities of certain toric singularities
can be generalized to arbitrary dimension, but we only need the case 
of dimension three.}
\end{rem}

Now we consider in detail the case when $n$ is a power of a prime number,
and the group $H_1$ is cyclic.

\begin{prop}
{ Let $K=K_{uvw}$ be a semigroup, defined by the conditions $\alpha u+
\beta v +\gamma w =0({\rm mod}{\em p^s})$ and $\alpha,\beta,\gamma\in 
\ZZ_{\geq 0},$ where $u$, $v$, and $w$ are some natural numbers.
The number $\delta$ defined in \ref{appdelta} equals
$p^{-s}{\rm min}_{K-\{0\}}(\alpha+\beta+\gamma).$ Then the number of homogeneous
triples $(u:v:w)$ such that $\delta(u,v,w)\geq \epsilon$ is at most
$2^2\epsilon^{-8}[4\epsilon^{-5}]_p$.}
\label{appfinmany}
\end{prop}

{\em Proof.} Consider the intersection of $K$ and the coordinate plane
$\alpha = 0$. It is the semigroup $K_1$ defined by the conditions
$\beta,\gamma\in \ZZ_{\geq 0},~\beta v+\gamma w = 0({\rm mod}{\em p^s}).$
If $\delta(u,v,w)\geq \epsilon$, then $\beta+\gamma\geq \epsilon p^s$
for all nonzero $(\beta,\gamma)\in K_1$. Therefore, the area of
$\RR^2_{>0}-conv(K_1-\{0\})$ is at least $\epsilon^2p^{2s}$, if
the area of the basic triangle in $\ZZ^2$ is equal to one. Because
any triangle in $\ZZ^2$ with no lattice points inside and on the
edges is basic, the number of points of $K_1$ that lie inside the
positive quadrant and on the boundary of $conv(K_1-\{0\})$ is at least
$-1+\epsilon^2p^{2s}/|\ZZ^2:span(K_1)|\geq -1+\epsilon^2p^s$.

The function $\beta-\gamma$ is monotone on the boundary of 
$conv(K_1-\{0\})$, and changes by at most $2p^s$ inside the positive
quadrant.  Hence, there is a segment of this boundary, that is represented
by the vector $(\beta_1,-\gamma_1)$ with $0<\beta_1,\gamma_1,~\beta_1
+\gamma_1\leq 2\epsilon^{-2}$. Hence there holds 
$v\beta_1=w\gamma_1({\rm mod}{\em p^s})$ with $0<\beta_1,\gamma_1,~\beta_1
+\gamma_1\leq 2\epsilon^{-2}$.

Analogously, we have $u\alpha_2=w\gamma_2({\rm mod}{\em p^s})$ with 
$0<\alpha_2,\gamma_2,~\alpha_2+\gamma_2\leq 2\epsilon^{-2}$.
Besides, ${\rm g.c.d.}(w,p^s)\leq[\epsilon^{-1}]_p$, because otherwise 
$(0,0,p^s/{\rm g.c.d.}(w,p^s))$ lies in $K$ and gives $\delta<\epsilon$.

There are at most $[\epsilon^{-1}]_p$ choices of $w({\rm mod}{\em p^s})$
up to multiplication by $(\ZZ/p^s\ZZ)^*$. There are at most
$2^2\epsilon^{-8}$ choices for the fourtuple $(\beta_1,\gamma_1,
\alpha_2,\gamma_2)$. Once we know $(w,\beta_1,\gamma_1,
\alpha_2,\gamma_2)$, there are at most $[2\epsilon^{-2}]_p$
for each of the numbers $u,v({\rm mod}{\em p^s})$. This proves the proposition.
\hfill$\Box$


\bigskip

\begin{thebibliography}{99}

\bibitem{Bolza} O. Bolza, {\em On binary sextics with linear transformations
into themselves}, Amer. J. Math., 1888, v.10, p.70.

\bibitem{CKM} H. Clemens, J. Koll\'ar, S. Mori, {\em Higher dimensional
 complex geometry}, Asterisque, 1988, v.166.

\bibitem{Hochster} D. Eisenbud, M. Hochster, {\em A Nullstellensatz with
 nilpotents and Zariski's main lemma on holomorphic functions},
 J. Algebra, 1979, v.58, N.1, pp.157-161.

\bibitem{Grady} K. O'Grady, {\em On the Kodaira dimension of moduli spaces
of abelian surfaces}, Compositio Math., 1989, v.72, N.2, pp.121-163.

\bibitem{Gritsenko} V. Gritsenko, {\em Modular forms and moduli spaces
of abelian and K3 surfaces}, Rossiiskaya Akademiya Nauk, Algebra i Analiz,
1994, v.6, N.6, pp.65-102.

\bibitem{GeerII} G.van der Geer, {\em Note on abelian schemes of level three},
Math. Ann., 1987, v.278, N.1-4, pp.401-408.

\bibitem{Geer} G. van der Geer, {\em On the geometry of a Siegel modular
 threefold},
 Math. Ann., 1982, Bd.260, N.3, pp.317-350.

\bibitem{Hammond} W.P. Hammond, {\em On the graded ring of Siegel modular
 forms of genus two}, Amer J. Math., 1965, v.87, N.2, pp.502-506.

\bibitem{Hulek} K. Hulek, G.K. Sankaran {\em The Kodaira dimension of certain
moduli spaces of abelian surfaces}, Compositio Math., 1994, v.90, N.1, pp.1-35.

\bibitem{Igusa} J.-I. Igusa, {\em A desingularization problem in the theory
 of Siegel modular functions}, Math. Ann., 1967, v.168, pp.228-260.

\bibitem{Lee} R. Lee, S. Weintraub, {\em An interesting algebraic variety},
Math. Intelligencer, 1986, v.8, N.1, pp.34-39.

\bibitem{Reid} M. Reid, {\em Canonical threefolds}, G\'eometrie Alg\'ebrique
 Angers. Sijthoff \&Noordhoff, 1980, pp.273-310.

\bibitem{Teissier} B. Teissier, {\em Mon\^omes, volumes et
multiplicit\'es},
Introduction \`a la th\'eorie des singularit\'es. II., Hermann, Paris,
1988, p.131.

\bibitem{Thompson} J.G. Thompson, {\em A finiteness theorem for subgroups of
 ${\rm PSl(2,\RR)}$ which are commensurable with ${\rm PSl(2,\ZZ)}$}, Proc. of
Symp.
 in Pure Math., 1980, v.37, pp.533-55.

\bibitem{Yamazaki} T. Yamazaki, {\em On Siegel modular forms of degree two},
 Amer J. Math., 1976, v.98, N.1, pp.39-53.

\end{thebibliography}
\end{document}